\input ./style/arxiv-general.cfg
\documentclass[aos,MSNbibl,seceqn,dvips]{arximspdf}
\makeatletter
   \@ifpackageloaded{graphicx}{}{\usepackage{graphicx}}
\makeatother
\usepackage{dcolumn}

\doi{10.1214/15-AOS1353}% Updated by VTEXPTS2LaTeX.exe, 25.08.2015
\volume{43}
\issue{6}
\pubyear{2015}
\firstpage{2588}
\lastpage{2623}
\docsubty{FLA}

\makeatletter
\newcolumntype{d}[1]{D{.}{.}{#1}}

\newcommand{\rrvert}{\vert}

\newcommand{\llvert}{\vert}
\newtheorem{thmm}{Theorem}[section]
\newtheorem{lem}[thmm]{Lemma}
\newtheorem{pro}[thmm]{Proposition}
\newtheorem{cor}[thmm]{Corollary}
\newproclaim{defi}[thmm]{Definition}
\newproclaim{rem}[thmm]{Remark}
\newproclaim{con}[thmm]{Condition}
\newproclaim{Exa}{Example}[section]
\makeatother

\begin{document}
\begin{frontmatter}

\title{Spectral statistics of large dimensional Spearman's rank
correlation matrix and its application}
\runtitle{Spectral statistics of Spearman's rank correlation matrix\hspace*{7pt}}

\begin{aug}
% Corresponding author: Wang Zhou - stazw@nus.edu.sg% Updated by
%VTEXPTS2LaTeX.exe, 25.08.2015 12:35
\author[A]{\fnms{Zhigang}~\snm{Bao}\thanksref{T1,m1}\ead[label=e1]{maomie2007@gmail.com}},
\author[B]{\fnms{Liang-Ching}~\snm{Lin}\thanksref{m2}\ead[label=e2]{lcln@mail.ncku.edu.tw}},
\author[A]{\fnms{Guangming}~\snm{Pan}\thanksref{T2,m1}\ead[label=e3]{gmpan@ntu.edu.sg}}
\and
\author[C]{\fnms{Wang}~\snm{Zhou}\corref{}\thanksref{T3,m3}\ead[label=e4]{stazw@nus.edu.sg}}
\runauthor{Bao, Lin, Pan and Zhou}
\affiliation{Nanyang Technological University\thanksmark{m1},
National Cheng Kung University\thanksmark{m2} and National University
of Singapore\thanksmark{m3}}
%\dedicated{}
\thankstext{T1}{Supported in part by the Ministry of Education,
Singapore, under Grant \# ARC 14/11, and NSF of China, Grant No. 11371317.}
\thankstext{T2}{Supported in part by a MOE Tier 2 Grant 2014-T2-2-060
and by a MOE Tier 1 Grant RG25/14 at the Nanyang Technological
University, Singapore.}
\thankstext{T3}{Supported in part by Grant R-155-000-151-112 at the
National University of Singapore.}

\address[A]{Z. Bao\\
G. Pan\\
Division of Mathematical Sciences\\
School of Physical and Mathematical Sciences\\
Nanyang Technological University\\
Singapore 637371\\
Singapore\\
\printead{e1}\\
\phantom{E-mail:\ }\printead*{e3}}

\address[B]{L.-C. Lin\\
Department of Statistics\\
National Cheng Kung University\\
Taiwan 70101\\
China\\
\printead{e2}}

\address[C]{W. Zhou\\
Department of Statistics and Applied Probability\\
National University of Singapore\\
Singapore 117546\\
Singapore\\
\printead{e4}}
\end{aug}

% HISTORY:
%
\received{\smonth{10} \syear{2014}}% Updated by VTEXPTS2LaTeX.exe,
%25.08.2015 12:35
%
\revised{\smonth{6} \syear{2015}}% Updated by VTEXPTS2LaTeX.exe,
%25.08.2015 12:35

% ABSTRACT
%
\begin{abstract}
Let $\mathbf{Q}=(Q_1,\ldots,Q_n)$ be a random vector drawn from the
uniform distribution on the set of all $n!$ permutations of $\{
1,2,\ldots,n\}$. Let $\mathbf{Z}=(Z_1,\ldots,Z_n)$, where $Z_j$ is the
mean zero variance one random variable obtained by centralizing and
normalizing $Q_j$, $j=1,\ldots,n$. Assume that $\mathbf
{X}_i,i=1,\ldots
,p$ are i.i.d. copies of $\frac{1}{\sqrt{p}}\mathbf{Z}$ and $X=X_{p,n}$
is the $p\times n$ random matrix with $\mathbf{X}_i$ as its $i$th row.
Then $S_n=XX^*$ is called the $p\times n$ Spearman's rank correlation
matrix which can be regarded as a high dimensional extension of the
classical nonparametric statistic Spearman's rank correlation
coefficient between two independent random variables. In this paper, we
establish a CLT for the linear spectral statistics of this
nonparametric random matrix model in the scenario of high dimension,
namely, $p=p(n)$ and $p/n\to c\in(0,\infty)$ as $n\to\infty$. We
propose a novel evaluation scheme to estimate the core quantity in
Anderson and Zeitouni's cumulant method in
[\textit{Ann. Statist.} \textbf{36} (2008) 2553--2576] to bypass the
so-called joint cumulant summability. In addition, we raise a {\emph
{two-step comparison approach}} to obtain the explicit formulae for the
mean and covariance functions in the CLT. Relying on this CLT, we then
construct a distribution-free statistic to test complete independence
for components of random vectors. Owing to the nonparametric property,
we can use this test on generally distributed random variables
including the heavy-tailed ones.
\end{abstract}

% KEYWORDS
\begin{keyword}[class=AMS]
\kwd[Primary ]{15B52}
\kwd[; secondary ]{62H10}
\end{keyword}
\begin{keyword}
\kwd{Spearman's rank correlation matrix}
\kwd{nonparametric method}
\kwd{linear spectral statistics}
\kwd{central limit theorem}
\kwd{independence test}
\end{keyword}
\end{frontmatter}

%s1 #&#
\section{Introduction}\label{sec1}
%s1.1 #&#
\subsection{Matrix model}\label{sec1.1}
In this paper, we will consider the large dimensional Spearman's rank
correlation matrices. First, we give the definition of the matrix
model. Let $P_n$ be the set consisting of all permutations of $\{
1,2,\ldots,n\}$. Suppose that $\mathbf{Z}=(Z_1,\ldots, Z_n)$ is a random
vector, where
\[
Z_i=\sqrt{\frac{12}{n^2-1}}\biggl(Q_i-\frac{n+1}{2}
\biggr)
\]
and $\mathbf{Q}:=(Q_1,\ldots, Q_n)$ is uniformly distributed on $P_n$.
That is, for any permutation $(\sigma(1),\sigma(2),\ldots, \sigma
(n))\in P_n$, one has
$
\mathbb{P}\{\mathbf{Q}=(\sigma(1),\sigma(2),\ldots, \sigma(n))\}=1/n!$.
For simplicity, we will use the notation $[N]:=\{1,\ldots,N\}$ for any
positive integer $N$ in the sequel. Now for $m\in[n]$ we conventionally
define the set of {\emph{$m$ partial permutations of $n$}} as
\[
P_{nm}:=\bigl\{(v_1,\ldots,v_m):
v_{1},\ldots,v_m\in[n]\mbox{ and } v_i\neq
v_j \mbox{ if }i\neq j \bigr\}.
\]
For any mutually distinct numbers $l_1,\ldots,l_m\in[n]$, it is
elementary to check that $(Q_{l_1},\ldots,Q_{l_m})$ is uniformly
distributed on $P_{nm}$. Such a fact immediately leads to the fact that
$\{Z_i\}_{i=1}^n$ is strictly stationary.
In addition, by setting $m=1$ or $2$, it is straightforward to check
through calculations that
%
%e1.1 #&#
\begin{equation}
\mathbb{E}Z_{i}=0,\qquad \mathbb{E}Z^2_{i}=1,\qquad
\operatorname {Cov}(Z_j,Z_k)=-\frac{1}{n-1}\qquad\mbox{if } j\neq k. \label{1113.20}
\end{equation}
Moreover, it is also easy to see that for any positive integer $l$,
\[
\mathbb{E}|Z_1|^{l}\leq C_l
\]
for some positive constant $C_l$ depending on $l$. Besides, we note
that $Z_i,i\in[n]$ are symmetric random variables.

Assuming that $\mathbf{X}_i=(x_{i1},\ldots, x_{in}),i=1,\ldots, p$ are
i.i.d. copies of $\frac{1}{\sqrt{p}}\mathbf{Z}$, we set $S_n=XX^*$, where
\[
X:=\left[ %
\matrix{\mathbf{X}_1
\cr
\vdots
\cr
\mathbf{X}_p }
 \right]=(x_{ij})_{p,n}.
\]
Then $S_n$ is referred to as the {\emph{Spearman's rank correlation matrix}}.

%s1.2 #&#
\subsection{Motivation from nonparametric statistics}\label{sec1.2}
A main motivation of considering the matrix $S_n$ is from
nonparametric statistics. We consider the hypothesis testing problem
on some random variable sequence $Y_1,\ldots, Y_p$ as
\begin{eqnarray*}
&&\mathbf{H}_0\!: \qquad Y_1,\ldots, Y_p \mbox{ are
independent}\quad \mbox{v.s.}\\
&&\mathbf{H}_1\!:\qquad Y_1,\ldots,
Y_p \mbox{ are not independent}.
\end{eqnarray*}
Note that since the covariance matrix and correlation matrix can
capture the dependence for Gaussian variables, it is
natural to compare the sample covariance or correlation matrix with
diagonal matrices to detect whether $\mathbf{H}_0$ holds in the
classical setting of {\emph{large $n$ and fixed $p$}}. Unfortunately,
due to the so-called {\emph{curse of dimensionality}}, it is now well
known that there is no hope to approximate the population covariance or
correlation matrix by sample ones in the situation of {\emph{large $n$
and comparably large $p$}} without any assumption imposed on the
population covariance matrix. However, it is still
possible to construct independence test statistics from the sample
covariance or correlation matrix for Gaussian variables even in the
high dimensional case. In the scenario of large $n$ and comparably
large $p$, there is a long list of literature devoted to studying the
properties of sample covariance matrix or sample correlation matrix
under the null hypothesis $\mathbf{H}_0$. % in varying degrees.
For example, Bao, Pan and Zhou \cite{BPZ2012}, as well as Pillai and
Yin \cite{PY20121}, studied the largest eigenvalue of sample
correlation matrix; Ledoit and Wolf \cite{LW2002} raised a quadratic
form of the spectrum of sample covariance matrix; Schott \cite
{Schott2005} considered the sum of squares of sample correlation
coefficients; Jiang \cite{Jiang2004} discussed the largest off-diagonal
entry of sample correlation matrix; Jiang and Yang \cite{JY2013}
studied the likelihood ratio test statistic for the sample correlation
matrix. However, for non-Gaussian variables, even in the classical
{\emph{large $n$ and fixed $p$}} case, the idea to compare the
population covariance matrix with diagonal matrices is substantially
invalid for independence test for those uncorrelated but dependent
variables. Moreover, for random vectors containing at least one heavy
tailed component such as a Cauchy random variable, there is even no
population covariance matrix. %For the high dimensional data, the
%philosophy of detecting dependence among variables via studying the
%statistics constructed from sample covariance or correlation matrix is
%also doubtful.

In view of the above, we discuss a nonparametric matrix model in this
paper and study its spectrum statistics under $\mathbf{H}_0$ in order
to % tackle the independence hypothesis test for
accommodate random variables with general distributions. Assume that
we have $n$ observations of the vector $(Y_1,\ldots, Y_p)$. Let
$Y_{11},\ldots, Y_{1n}$ be the observations of the first coefficient
$Y_1$ and set $Q_{1j}$ to be the rank of $Y_{1j}$ among $Y_{11},\ldots,
Y_{1n}$. We then replace
$(Y_{11},\ldots, Y_{1n})$ by the corresponding normalized rank sequence
$(x_{11},\ldots,x_{1n})$, where
\[
x_{1j}=\sqrt{\frac{12}{p(n^2-1)}}\biggl(Q_{1j}-
\frac{n+1}{2}\biggr),\qquad j\in[n].
\]
Analogously, we can define the rank sequence $(x_{i1},\ldots, x_{in})$
for other $i\in[p]$.
For simplicity, in this paper we only consider the case where $Y_i,i\in
[p]$ are continuous random variables. In this case, the probability of
a tie occurring in the sequence $Q_{i1},\ldots, Q_{in}$ for any $i\in
[p]$ is zero. Then $S_n=XX^*$ with $X=(x_{ij})_{p,n}$ is the so-called
Spearman's rank correlation matrix under $\mathbf{H}_0$, which can be
regarded as a high dimensional extension of the classical Spearman's
rank correlation coefficient between two random variables. Then we can
construct statistics from $S_n$ to tackle the above hypothesis testing
problem. By contrast, the parametric models such as Pearson's sample
correlation matrix and sample covariance matrix are well studied by
statisticians and probabilists. However, the work on Spearman's rank
correlation matrix is few and far between. In \cite{BZ2008}, Bai and
Zhou proved that the limiting spectral distribution of $S_n$ is also
the famous Marchenko--Pastur law (MP law). In~\cite{Zhou2007}, Zhou
studied the limiting behavior of the largest off-diagonal entry of
$S_n$. Our purpose in this paper is to derive the fluctuation (a CLT)
of the linear spectral statistics for $S_n$. As an application, we will
construct a nonparametric statistic to detect dependence of components
of random vectors.
%s1.3 #&#
\subsection{Main result}\label{sec1.3}
We set $\lambda_1\geq\cdots\geq\lambda_p$ to be the ordered
eigenvalues of $S_n$. Our main task in this paper is to study the
limiting behavior of the so-called linear spectral statistics
$
\mathcal{L}_n[f]=\sum_{i=1}^pf(\lambda_i)
$
for some test function $f$. In this paper, we will focus on the
polynomial test functions and, therefore, it suffices to study the
joint limiting behavior of $\operatorname{tr} S^{k}_n, k=1,\ldots
,\infty$. We
state the main result as the following theorem.

%th1.1 #&#
\begin{thmm} \label{thmm.1111.1}Assuming that both $n$ and $p:=p(n)$ tend
to $\infty$ and
\[
n/p\rightarrow c\in(0,\infty),
\]
we have
\[
\bigl\{\operatorname{tr} S^{k}_n-\mathbb{E}
\operatorname{tr} S^{k}_n \bigr\} _{k=2}^\infty
\Longrightarrow\{G_k\}_{k=2}^\infty\qquad \mbox{as } n
\rightarrow\infty,
\]
where $\{G_k\}_{k=2}^\infty$ is a Gaussian process with mean zero and
the covariance function given by
%
%e1.2 #&#
\begin{eqnarray}
\label{1111.2}&& \operatorname{Cov}(G_{k_1},G_{k_2})\nonumber\\
&&\qquad=2c^{k_1+k_2}
\sum_{j_1=0}^{k_1-1}\sum
_{j_2=0}^{k_2}\pmatrix{k_1
\cr
_1}\pmatrix{k_2
\cr
j_2} \biggl(
\frac{1-c}{c} \biggr)^{j_1+j_2}
\nonumber
\\
&&\qquad\quad{} \times\sum_{l=1}^{k_1-j_1}l\pmatrix
{2k_1-1-(j_1+l)
\cr
k_1-1}
\pmatrix{2k_2-1-j_2+l
\cr
k_2-1}
\\
&&\qquad\quad{} -2c^{k_1+k_2+1}\sum_{j_1=0}^{k_1}\sum
_{j_2=0}^{k_2}\pmatrix {k_1
\cr
j_1}\pmatrix{k_2
\cr
j_2} \biggl(
\frac{1-c}{c} \biggr)^{j_1+j_2}\nonumber\\
&&\qquad\quad{}\times\pmatrix {2k_1-j_1
\cr
k_1-1}\pmatrix{2k_2-j_2
\cr\nonumber
k_2-1}.
\end{eqnarray}
Moreover, we have the following expansion for the expectation function:
%
%e1.3 #&#
\begin{eqnarray}
\label{1111.3} \mathbb{E}\operatorname{tr} S^{k}_n&=&
\frac{n^k}{(n-1)^{k-1}}\sum_{j=0}^{k-1}
\frac{1}{j+1}\pmatrix{k
\cr
j}\pmatrix{k-1
\cr
j} \biggl(\frac
{n-1}{p}
\biggr)^{j}-\frac{1}{2}\sum_{j=0}^k
\pmatrix{k
\cr
j}^2c^j
\nonumber
\\
&&{}+2c^{1+k}\sum_{j=0}^k
\pmatrix{k
\cr
j} \biggl(\frac{1-c}{c} \biggr)^j \biggl[\pmatrix{2k-j
\cr
k-1}-\pmatrix{2k+1-j
\cr
k-1} \biggr]
\\
&&{}+\frac{1}{4} \bigl[(1-\sqrt{c})^{2k}+(1+\sqrt{c})^{2k}
\bigr]+o(1).\nonumber
\end{eqnarray}
\end{thmm}

%re1.2 #&#
\begin{rem} Note that when $k=1$, $\operatorname{tr} S_n=n$ is deterministic.
Actually, it can also be checked that the right-hand side of (\ref
{1111.2}) is zero when $k_1=k_2=1$ and the right-hand side of (\ref
{1111.3}) equals $n+o(1)$ when $k=1$ therein.
\end{rem}

%s1.4 #&#
\subsection{Methodologies of the proof}\label{sec1.4}
Roughly speaking, we will start with a {\emph{cumulant method}}
introduced by Anderson and Zeitouni in \cite{AZ2009} to establish the
CLT for $\{(\operatorname{tr} S^{k}_n-\mathbb{E}\operatorname{tr}
S^{k}_n)/\sqrt{\operatorname
{Var}(\operatorname{tr} S^{k}_n})\}_{k=2}^{\infty}$. The cumulant
method can be
viewed as a modification of the celebrated moment method in Random
Matrix Theory (RMT). Without trying to be comprehensive, we refer to
\cite{BS2004,BS2008,SS1998,Wigner1955} for further reading.
Particularly, we refer to \cite{AZ2009,CL1995,Soshnikov2002} for the
cumulant method in proving CLTs for linear spectral statistics in RMT.

As explained at the beginning, $\{Z_i\}_{i=1}^n$ is a stationary
sequence, which inspires us to learn from \cite{AZ2009}. However, in
\cite{AZ2009} the stationary sequence is required to satisfy the
so-called {\emph{joint cumulant summability}} property [see (\ref
{1111.4}) below], which was used to bound the core quantity (\ref
{1012.3}). The property of joint cumulant summability is crucial in the
spectral analysis of time series; one can refer to \cite
{Brillinger2001,Rosenblatt1985,WS2004,XW2011} for further reading.
Unfortunately, to verify this property for a general stationary
sequence is highly nontrivial. In this paper, we will not try to check
whether the joint cumulant summability holds for $\mathbf{Z}$. Instead
we will provide a relatively rough but crucial bound on (\ref{1012.3})
through a totally novel evaluation scheme; see Proposition~\ref
{lem.106.15} and Corollary~\ref{cor.107.20} below. Such a bound will
allow us to bypass the joint cumulant summability and serve as a main
input to pursue the cumulant method to establish the CLT.
With this CLT, what remains is therefore to evaluate the nonnegligible
terms of the mean and covariance functions.
It will be shown that the mean and covariance functions of $\{
\operatorname{tr}
S^{k}_n\}_{k=2}^{\infty}$ can be expressed by some sums of terms
indexed by set partitions. For the explicit values of the
nonnegligible terms of the mean and covariance functions, we will
adopt a {\emph{two-step comparison strategy}} to compare these
expressions with the existing results for the sample covariance matrices.

%s1.5 #&#
\subsection{Organization and notation}\label{sec1.5}
Our paper is organized as follows. Section~\ref{sec2} is devoted to the
application of our CLT to the independence test on random vectors. In
Section~\ref{sec3}, we will introduce some basic notions of joint
cumulants and
some known results from \cite{AZ2009}. Section~\ref{sec4} is our main
technical
part which will be devoted to providing the required bound for the sum
(\ref{1012.3}). Specifically, we will mainly prove Proposition~\ref
{lem.106.15} therein. In Section~\ref{sec5}, we will use the bounds
obtained in
Section~\ref{sec4} to show that all high order cumulants tend to zero
when $n\to
\infty$. Finally, in Section~\ref{sec6}, we will combine the bounds
in Section~\ref{sec4}
and the aforementioned two-step comparison strategy to evaluate the
main terms of the mean and covariance functions.

Throughout the paper, we will use $\#\mathbb{S}$ to represent the
cardinality of a set $\mathbb{S}$. For any number set $\mathbb
{A}\subset[n]$, we will use $\{j_\alpha\}_{\alpha\in\mathbb{A}}$ to
denote the set of $j_\alpha$ with $\alpha\in\mathbb{A}$. Analogously,
$(j_\alpha)_{\alpha\in\mathbb{A}}$ represents the vector obtained by
deleting the components $j_\beta$ with $\beta\in[n]\setminus\mathbb
{A}$ from $(j_1,\ldots, j_n)$. In addition, we will use the notation
$\mathrm{i}=\sqrt{-1}$ to denote the imaginary unit to release $i$
which will be frequently used as subscript or index. For any vector
$\vec{\xi}=(\xi_1,\ldots,\xi_N)$, we say the {\emph{position}} of
$\xi
_i$ in $\vec{\xi}$ is $i$. For example, for vector
$(Z_{j_2},Z_{j_1},Z_{j_4},Z_{j_3})$, the position of $Z_{j_1}$ is $2$.
Moreover, we will use $C$ to represent some positive constant
independent of $n$ whose value may be different from line to line.

%s2 #&#
\section{Application on independence test}\label{sec2}

In this section, we consider an application of Theorem~\ref{thmm.1111.1}.
We construct a nonparametric statistic to test complete independence
for the components of a high dimensional random vector.
Our proposed statistic is $W_7$ and we highlight it here at first
\begin{eqnarray*}
W_7&=&W_7(k,\delta)\\
&:=&\frac{\operatorname{tr} S_n^k-\mathbb
{E}\operatorname{tr}
S_n^k}{\sqrt{\operatorname{Var}(G_k)}}+
n^{-\delta} \biggl[n \biggl(\max_{1\leq i < j \leq p} \biggl\llvert
\frac
{p}{n}s_{ij} \biggr\rrvert \biggr)^2 -4\log p+
\log\log p \biggr],
\end{eqnarray*}
where $ 0 < \delta< 1$. Note that our statistic depends on the
parameters $k$ and $\delta$. We will discuss how to choose $k$ and
$\delta$ at the end of this section. To see the performance of our
statistic, we will do
the simulation under various settings.
We also compare our statistic with some other parametric or
nonparametric statistics in the literature.
To introduce these statistics, we shall define some notation at first.
Let $R_n = (r_{ij})_{p\times p}$ be the Pearson's sample correlation
matrix based on $n$ independent copies of
$p$-dimensional random vector $(Y_1,\ldots,Y_p)$. We denote by
$\lambda
_{\operatorname{max}}(R)$ the largest eigenvalue of $R_n$.
In addition, we denote by $s_{ij}$ the $(i,j)$th entry of $S_n$.
We compare the performance of our statistic $W_7$ with the following statistics:
\begin{eqnarray}
\mathrm{(i)}\quad W_1& =&\frac{n\lambda_{\operatorname
{max}}(R)-(p^{1/2}+n^{1/2})^2}{(n^{1/2}+p^{1/2})(p^{-1/2}+n^{-1/2})^{1/3}}
\qquad \mbox{(see \cite{BPZ2012} or \cite{PY20121})},
\nonumber\\
\mathrm{(ii)}\quad  W_2&=& \frac{\operatorname{tr}
S_n^k-\mathbb{E}\operatorname{tr} S_n^k}{\sqrt{\operatorname
{Var}(G_k)}}\qquad \mbox{(see Theorem~\ref{thmm.1111.1})},
\nonumber\\
\mathrm{(iii)}\quad W_3& =& \frac{\sum_{i=2}^{p}\sum_{j=1}^{i-1} r_{ij}^2 - p(p-1)/(2n)}{(p/n)}\qquad \mbox{(see \cite{Schott2005})},
\nonumber\\
\mathrm{(iv)}\quad W_4& =& \frac{\log(|R_n|) -
(p-n+3/2)\log(1-{p}/{(n-1)})+(n-2({p}/{(n-1)}))}{\sqrt{-2[
{p}/{(n-1)} + \log(1-{p}/{(n-1)})]}}\nonumber\\
 \eqntext{\mbox{(see \cite{JY2013})},}
\\
\mathrm{(v)}\quad  W_5& =& n \biggl(\max_{1\leq i < j
\leq p} \llvert  r_{ij} \rrvert  \biggr)^2 -4\log n+ \log\log n\qquad \mbox{(see
\cite
{Jiang2004})},\nonumber
\\
\mathrm{(vi)}\quad W_6 &=& n \biggl(\max_{1\leq i <
j \leq p}\biggl \vert  \frac{p}{n}s_{ij} \biggr\vert  \biggr)^2 -4\log p+ \log
\log
p\qquad\mbox{(see \cite{Zhou2007})}.\nonumber
\end{eqnarray}

At first, from Theorem~\ref{thmm.1111.1}, we see that the CLT holds for
the statistic $W_2$.
By our construction, $W_7 = W_2 + n^{-\delta} W_6$, which can be
regarded as a slight modification of $W_2$ by adding a small penalty in
terms of $W_6$.
We expect that the statistic $W_7$ will take the advantages of both
$W_2$ and $W_6$,
and thus its performance will be better. More specifically, we
illustrate the philosophy of such a construction via the following
examples. Let $\mathbf{g}$ be a $p$-dimensional Gaussian vector with
the population covariance matrix $\Sigma_\mathbf{g}$.
Two extreme alternative hypotheses are considered below.
The first case is that $\Sigma_\mathbf{g}$ has only one significantly
large off-diagonal entry.
Then the corresponding Spearman's rank correlation matrix will also
have a significantly large off-diagonal entry.
Since $W_6$ is constructed from the largest off-diagonal entry,
it is sensitive to this kind of dependence structure.
In contrast, one cannot tell the dependence structure in $\mathbf{g}$
by $W_2$ since
the linear spectral statistics are relatively robust under the
disturbance of a single entry of the population covariance matrix.
However, $W_7$ has an additional penalty which is $n^{-\delta} W_6$
compared to $W_2$. So one can capture the dependence contained in
$\mathbf{g}$ by $W_7$.
Now, we consider the second case where $\Sigma_\mathbf{g}$ contains a
lot of small nonzero off-diagonal entries.
In this case, the statistic $W_6$ performs badly since the largest
off-diagonal entry of $\Sigma_\mathbf{g}$ is close to zero.
In contrast, the statistic $W_7$ performs as well as $W_2$ since the
spectral statistics can accumulate all the effects caused by
these small off-diagonal entries.
%We remark here that in the high dimensional case one cannot regard
%$\Sigma_\mathbf{g}$ to be approximately equal to identity if there are
%a lot of small off-diagonal entries.

We below summarize the limiting null distributions of $W_i$,
$i=1,\ldots
,7$, and the corresponding assumptions in the references \cite
{BPZ2012,PY20121,Schott2005,JY2013,Jiang2004,Zhou2007}.
The null distribution of $W_1$ converges to the type 1 Tracy--Widom
law (see \cite{BPZ2012,PY20121}), assuming that the variables
$Y_1,\ldots, Y_p$ have sub-exponential tails.
The limiting null distribution of $W_2$ is $N(0,1)$ by Theorem~\ref
{thmm.1111.1}.
In \cite{Schott2005} and \cite{JY2013}, the weak convergence of $W_3$
and $W_4$ to $N(0,1)$ is established for the Gaussian vector
$(Y_1,\ldots, Y_p)$ only.
If one assumes that $Y_1,\ldots, Y_p$ are i.i.d. with
$\mathbb{E}Y_1^{30-\varepsilon}<\infty$ for any constant
$\varepsilon>0$,
the limiting distribution of $W_5$ derived in \cite{Jiang2004}
possesses the following c.d.f.:
\[
F_{W_5}(y) = \exp \bigl\{ -\bigl(c^2 \sqrt{8\pi}
\bigr)^{-1} e^{-y/2} \bigr\},
\]
which is called the extreme distribution of type I. According to \cite
{Zhou2007}, the statistic $W_6$ is distribution-free,
with the following asymptotic distribution:
\[
F_{W_6}(y) = \exp \bigl\{ -(8\pi)^{-1/2} e^{-y/2}
\bigr\}.
\]
Clearly, the statistic $W_7 = W_2 + n^{-\delta} W_6$ is also
distribution-free and possesses the limiting distribution
$N(0,1)$ due to Slutsky's theorem. %since $n^{-\delta} \rightarrow0$

We denote by $\operatorname{Cauchy}(\alpha,\beta)$ the Cauchy distribution
with location parameter $\alpha$ and scale parameter $\beta$.
In addition, we denote by $t(\gamma)$ the student's $t$-distribution
with degrees of freedom $\gamma$.
We consider three null hypotheses with the nominal significance level
$\alpha= 5\%$:
\begin{itemize}
\item$\mathbf{H}_{0,1}$: $\mathbf{Y}_{j},j\in[n]$ are i.i.d.
$N_p(\mathbf{0},I_p)$ vectors;
\item$\mathbf{H}_{0,2}$: $Y_{ij},i\in[p],j\in[n]$ are i.i.d.
$\operatorname
{Cauchy}(0,1)$ variables.;
\item$\mathbf{H}_{0,3}$: $Y_{i_1,j}$ are i.i.d. $N(0,1)$ variables;
$Y_{i_2,j}$ are i.i.d.
$\operatorname{Cauchy}(0,1)$ variables; $Y_{i_3,j}$ are i.i.d. $t(4)$
variables,
where $i_1 = 1,\ldots,\lfloor p/3 \rfloor$, $i_2 = \lfloor p/3
\rfloor
+ 1,\ldots,\lfloor2p/3 \rfloor$,
$i_3 = \lfloor2p/3 \rfloor+ 1,\ldots,p$ and $j\in[n]$.
\end{itemize}
For each null hypothesis, we consider two alternatives:
\begin{itemize}
\item$\mathbf{H}_{a,1-1}$ (one large disturbance): $\mathbf
{Y}_{j},j\in
[n]$ are i.i.d. $N_p(\mathbf{0},I_p+C)$ and $C = (c_{ik})_{p\times p}$
with $c_{ik}=0$, $i,k\in[p]$, except $c_{12} = c_{21} =0.8$.
\item$\mathbf{H}_{a,1-2}$ (many small disturbances): $\mathbf
{Y}_{j},j\in[n]$ are i.i.d. $N_p(\mathbf{0},I_p+D)$ and $D =
(d_{ik})_{p\times p}$ with $d_{ik}=4/p$ except $d_{ii} = 0$, $i,k\in[p]$.
\item$\mathbf{H}_{a,2-1}$ (one large disturbance): $X_{ij}$ are i.i.d.
$\operatorname{Cauchy}(0,1)$. We set the observations $Y_{1j} = X_{1j} +0.8
X_{2j}$, $Y_{2j} = X_{2j} +0.8 X_{1j}$ and $Y_{ij} = X_{ij}$, for all
$i=3,\ldots,p$ and $j\in[n]$;
\item$\mathbf{H}_{a,2-2}$ (many small disturbances): $X_{ij}$ are
i.i.d. $\operatorname{Cauchy}(0,1)$. We set the observations $Y_{ij} =
X_{ij} +
(7p)^{-1}\sum_{k\neq i} X_{kj}$ for $i\in[p]$ and $j\in[n]$;
\item$\mathbf{H}_{a,3-1}$ (one large disturbance): the vectors
$(Y_{1,j},\ldots,Y_{{\lfloor p/3\rfloor},j}),j\in[n]$ are i.i.d.
$N_{\lfloor p/3\rfloor}(\mathbf{0},I_{\lfloor p/3\rfloor}+C')$ and $C'
= (c_{ik})_{\lfloor p/3\rfloor\times\lfloor p/3\rfloor}$ with
$c_{ik}=0$, $i,k=1,\ldots,\lfloor p/3\rfloor$, except $c_{12} = c_{21}
=0.8$. Moreover, $Y_{i,j},i=\lfloor p/3\rfloor+1,\ldots,n, j\in[n]$
are the same as those in $\mathbf{H}_{0,3}$.
\item$\mathbf{H}_{a,3-2}$ (many small disturbances): the vectors
$(Y_{1,j},\ldots,Y_{{\lfloor p/3\rfloor},j}),j\in[n]$ are i.i.d.
$N_{\lfloor p/3\rfloor}(\mathbf{0},I_{\lfloor p/3\rfloor}+D')$ and $D'
= (d_{ik})_{\lfloor p/3\rfloor\times\lfloor p/3\rfloor}$ with
$d_{ik}=12/p$ except $d_{ii} = 0$, $i,k=1,\ldots,\lfloor p/3\rfloor$.
Moreover, $Y_{ij},i=\lfloor p/3\rfloor+1,\ldots,n, j\in[n]$ are the
same as those in $\mathbf{H}_{0,3}$.
\end{itemize}

The results of sizes and powers listed in Table~\ref{results} are
based on the choices of
$(n,p) = (60,40)$, $(120,80)$, $(60,10)$ and $(120,160)$ and 1000 replications.
The tuning parameters of $W_2$ and $W_7$ are set to be $k=4$ and
$\delta=0.5$, which will be explained later.
In the case of $(n,p)=(120,160)$, $W_4$ is not defined. Hence, we
ignore it in Table~\ref{results}.
Moreover, in the cases of $\mathbf{H}_{0,2}$ and $\mathbf{H}_{0,3}$,
the distribution assumption or the moment assumption is violated for
the statistics $W_1$, $W_3$, $W_4$ and $W_5$.
Therefore, the corresponding values are also absent. We summarize our
findings as follows.
\begin{longlist}[(1)]
\item[(1)] The sizes of $W_2$, $W_3$, $W_4$ and $W_7$ are close to the
nominal size $5\%$.
However, $W_3$ has some size distortion in the case of $(n,p)=(120,160)$.
Meanwhile, the sizes of $W_1$, $W_5$ and $W_6$ tend to be smaller than
$5\%$.
\item[(2)] If the alternative hypothesis is the case of one large disturbance,
$W_5$, $W_6$ and $W_7$ outperform the other statistics.
In contrast, if the alternative hypothesis is the case of many small
disturbances,
$W_1$, $W_2$, $W_3$, $W_4$ and $W_7$ have better performance than
$W_5$ and $W_6$.
\end{longlist}

Overall, the size of $W_7$ is close to the nominal level $\alpha=5\%$
in our simulation study.
Moreover, $W_7$ has higher powers than the other statistics in most
cases of the alternative hypotheses.

%
%t1 #&#
\begin{table}
\tabcolsep=0pt
\caption{The sizes and powers (percentage) of $W_1$ to $W_7$ for
different hypotheses, sample size $n$ and dimension $p$}\label{results}
%\begin{center}
%
{\fontsize{8.6}{10.6}{\selectfont
\begin{tabular*}{\textwidth}{@{\extracolsep{\fill}}ld{3.1}d{3.1}d{3.1}d{3.1}d{3.1}d{3.1}d{3.1}d{3.1}d{3.1}d{3.1}d{3.1}d{3.1}d{3.1}@{}}
\hline
\multicolumn{1}{@{}l}{$\bolds{(n,p)}$} &
\multicolumn{1}{c}{$\bolds{W_1}$}&\multicolumn{1}{c}{$\bolds{W_2}$}&\multicolumn{1}{c}{$\bolds{W_3}$}&
\multicolumn{1}{c}{$\bolds{W_4}$}&\multicolumn{1}{c}{$\bolds{W_5}$}&\multicolumn{1}{c}{$\bolds{W_6}$}&\multicolumn{1}{c}{$\bolds{W_7}$}
 &\multicolumn{1}{c}{$\bolds{W_2}$}&\multicolumn{1}{c}{$\bolds{W_6}$}&\multicolumn{1}{c}{$\bolds{W_7}$}&\multicolumn{1}{c}{$\bolds{W_2}$}&
 \multicolumn{1}{c}{$\bolds{W_6}$}&\multicolumn{1}{c@{}}{$\bolds{W_7}$}\\
\hline
&\multicolumn{7}{c}{$\mathbf{H}_{0,1}$} &
\multicolumn{3}{c}{$\mathbf{H}_{0,2}$} &
\multicolumn{3}{c@{}}{$\mathbf{H}_{0,3}$}
\\[-6pt]
&\multicolumn{7}{c}{\hrulefill} &
\multicolumn{3}{c}{\hrulefill} &
\multicolumn{3}{c@{}}{\hrulefill}
\\
$(60,40)$ &0.4&4&5&4.9&2.1&3.2&4.5 &4.7&1.9&4.9&4.5&2.0&4.8\\
$(120,80)$ &1.4&5&5.6&4.9&2.6&2.2&5.7&4.6&2.6&5.4&3.7&3.7&4.4\\
$(60,10)$ &0.4&3.6&2.1&3.1&3.8&3.4&4.1&3.6&2.2&4.2&3.9&3.6&4.2\\
$(120,160)$ &1.7&5.8&10.9&\mbox{--}&1.9&2.3&6.7 &4&2.8&4.8&4.2&2.4&4.6\\[3pt]
&\multicolumn{7}{c}{$\mathbf{H}_{a,1-1}$}& \multicolumn {3}{c}{$\mathbf{H}_{a,2-1}$} & \multicolumn{3}{c}{$\mathbf
{H}_{a,3-1}$} \\[-6pt]
&\multicolumn{7}{c}{\hrulefill} &
\multicolumn{3}{c}{\hrulefill} &
\multicolumn{3}{c@{}}{\hrulefill}\\
%$(n,p)$ &$W_1$&$W_2$&$W_3$&$W_4$&$W_5$&$W_6$&$W_7$&
%&$W_2$&$W_6$&$W_7$&$W_2$&$W_6$&$W_7$&\\ \hline
$(60,40)$ &2.5&13.9&22.6&17.4&100&100&92.9
&25.6&100&99.7&12.4&99.9&98.2\\
$(120,80)$ &5.6&13.2&25&20.5&100&100&97.5
&26.2&100&100&12.1&100&96\\
$(60,10)$ &29.4&84.3&96.2&99.5&100&100&99.3
&100&100&100&83.4&100&99.5\\
$(120,160)$ &2.2&7.7&21.5&\mbox{--}&100&100&99.6
&11.4&100&100&6.5&100&94.6\\ [3pt]
&\multicolumn{7}{c}{$\mathbf{H}_{a,1-2}$}& \multicolumn
{3}{c}{$\mathbf{H}_{a,2-2}$} & \multicolumn{3}{c}{$\mathbf
{H}_{a,3-2}$}\\[-6pt]
&\multicolumn{7}{c}{\hrulefill} &
\multicolumn{3}{c}{\hrulefill} &
\multicolumn{3}{c@{}}{\hrulefill}\\ %\cline{3-9} \cline{12-14} \cline{17-19}
%$(n,p)$ &$W_1$&$W_2$&$W_3$&$W_4$&$W_5$&$W_6$&$W_7$&
%&$W_2$&$W_6$&$W_7$&$W_2$&$W_6$&$W_7$&\\ \hline
$(60,40)$ &99.9&99.8&99.7&67.5&15.6&18&99.8
&98.2&28.7&97.8&99.9&65.4&99.9\\
$(120,80)$ &100&100&100&72.4&11.8&12.3&100
&100&68.7&100&100&39&100\\
$(60,10)$ &100&100&100&100&99.8&99.5&100
&91.7&45.6&91.4&100&100&100\\
$(120,160)$ &100&99.3&99.2&\mbox{--}&4&4.3&98.9
&100&73.2&100&99.6&6.1&99.4\\
\hline
\end{tabular*}}}
%
%\end{center}
\end{table}

Finally, we consider how to choose the parameters $k$ and $\delta$ in
$W_2$ and $W_7$.
For illustration, we consider the case $\mathbf{H}_{0,1}$ versus
$\mathbf{H}_{a,1-1}$.
The parameter $k$'s are chosen to be $2,4,6,8$ and $10$. The parameter
$\delta$'s are chosen to be $0.3,0.4,0.5,0.6,0.7$ and $0.8$,
The sample size $n$ and the dimension $p$ are set to be $(60,40),
(120,80), (60,10)$ and $(120,160)$.
The sizes and powers of $W_2$ and $W_7$ are given in Tables~\ref
{sensize} and \ref{senpower}, respectively.
Based on Tables~\ref{sensize} and \ref{senpower}, for a fixed value of
$\delta$, one can see that
their sizes and powers are robust when $k\geq4$.
Therefore, we suggest to set $k$ as $4$ for both $W_2$ and $W_7$.
(Of course, theoretically, any $k$ larger than 4 is applicable, but it
will increase the computational burden without significant benefit).

The parameter $\delta$ should be chosen appropriately, so that $W_7$
will inherit the desirable properties of $W_2$ and $W_6$. In principle,
$\delta$ should not be very small or very large. If it is too close to
zero, that is, $n^{-\delta}$ is close to one, the size of $W_7$ will be
influenced since the limiting null distribution of $W_7$ may not be
standard normal any more. On the other hand, if $\delta$ is relatively
large, $W_7$ will not detect the dependent case where $W_6$ works. In
other words, the power of $W_7$ is weak.
In Tables~\ref{sensize} and \ref{senpower}, for $k=4$,
we can see that $W_7$ performs the best when $\delta$ is $0.5$ for
these four combinations of $(n,p)$.
We also conduct simulations under the other two null hypotheses
$\mathbf{H}_{0,2}$ and $\mathbf{H}_{0,3}$. The results are similar to
$\mathbf{H}_{0,1}$.
So based on our simulations we suggest to \textit{use $k=4, \delta
=0.5$}. This empirical choice is independent of the specific
distribution of $(Y_1,\ldots,Y_p)$ since both $W_2$ and $W_7$ are
distribution free.
%Thus the chosen of $\delta$ should not depend on the distribution of
%the original variables $Y_1,\ldots,Y_p$ under the null hypothesis.
% However, there are many combinations of $(n,p)$ and it may not have a
%rule or criterion to choose the optimal $\delta$ theoretically.
% We suggest the following way to choose the parameter $\delta$ in
%application.
% For a given combination of $(n,p)$, we can simulate the independent
%Gaussian variables
% to obtain the sizes of $W_7$ based on different values of $\delta$.
% Then, we can select the smallest $\delta$ based on the suitable range
%of the simulated sizes,
% say range from $0.4\sim0.6$ for nominal significant level $\alpha= 5
%\%$.

%
%t2 #&#
\begin{table}
\caption{The sizes (percentage) of $W_2$ and $W_7$ for different $n$,
$p$, $k$ and $\delta$ under the null hypothesis $\mathbf{H}_{0,1}$}\label{sensize}
\begin{tabular*}{\textwidth}{@{\extracolsep{\fill}}ld{2.0}d{1.1}cd{2.1}d{2.1}d{2.1}d{2.1}d{2.1}d{2.1}@{}}
\hline
&&& \multicolumn{7}{c@{}}{$\bolds{W_7}$}\\ [-6pt]
&&& \multicolumn{7}{c@{}}{\hrulefill}  \\
\multicolumn{1}{@{}l}{$\bolds{(n,p)}$} & \multicolumn{1}{c}{$\bolds{k}$} & \multicolumn{1}{c}{$\bolds{W_2}$} &
\multicolumn{1}{c}{$\bolds{\delta= }$} & \multicolumn{1}{c}{$\bolds{0.3}$} & \multicolumn{1}{c}{$\bolds{0.4}$} & \multicolumn{1}{c}{$\bolds{0.5}$} & \multicolumn{1}{c}{$\bolds{0.6}$} &
\multicolumn{1}{c}{$\bolds{0.7}$} & \multicolumn{1}{c}{$\bolds{0.8}$} \\
 \hline
$(60,40)$& 2 & 4.8 && 22.7 & 14.7 & 10.5 & 7.9 & 7 & 6.4  \\
& 4 & 4.5 && 19 & 12.7 & 6 & 5.2 & 4.9 & 4.9 \\
& $6$ & 4 && 16.4 & 9.6 & 5.8 & 5.4 & 4.7 & 4.2 \\
& 8 & 3.3 && 16.9 & 9 & 4.8 & 4.2 & 3.7 & 3.3 \\
& 10& 4 && 16.2 & 6.5 & 4 & 3.6 & 3.8 & 4  \\ [3pt]
$(120,80)$& 2 & 5.4 && 14.7 & 9.5 & 7.5 & 6.2 & 5.7 & 5.6  \\
& 4 & 5.2 && 14.3 & 8.9 & 5.5 & 5.4 & 5.3 & 5.3  \\
& $6$ & 4.1 && 12.7 & 6.3 & 4.8 & 4.4 & 4.1 & 4.1  \\
& 8 & 4.2 && 12 & 6.7 & 5.5 & 4.8 & 4.8 & 4.5  \\
& 10 & 3.8 && 13.3 & 7.2 & 5 & 4.2 & 4 & 4 \\ [3pt]
$(60,10)$& 2 & 3.9 && 21.3 & 13 & 5.9 & 5 & 4.8 & 4.3  \\
& 4 & 4.3 && 18.2 & 8.8 & 4.9 & 4.7 & 4.6 & 4.6 \\
& $6$ & 4.4 && 17.6 & 7.6 & 5.7 & 4.7 & 4.6 & 4.5 \\
& 8 & 4.9 && 15.5 & 7 & 5.6 & 5.3 & 5.1 & 5  \\
& 10 & 5 && 13.6 & 6.5 & 5 & 5.2 & 5.1 & 5.1  \\[3pt]
$(120,160)$& 2 & 6 && 22.7 & 12.4 & 9.3 & 8.2 & 7.2 & 6.7  \\
& 4 & 5.4 && 13.9 & 8 & 5.8 & 5.6 & 5.5 & 5.4  \\
& $6$ & 5.3 && 15.2 & 9.1 & 5.7 & 5.4 & 5.3 & 5.3  \\
& 8 & 5.3 && 13.4 & 8.6 & 5.8 & 5.6 & 5.6 & 5.4  \\
& 10 & 4.3 && 13.4 & 7.4 & 5.2 & 4.8 & 4.5 & 4.3  \\
\hline
\end{tabular*}
\end{table}
%

%
%t3 #&#
\begin{table}
\caption{The powers (percentage) of $W_2$ and $W_7$ for different $n$,
$p$, $k$ and $\delta$ under the alternative hypothesis~$\mathbf{H}_{a,1-1}$}\label{senpower}
\begin{tabular*}{\textwidth}{@{\extracolsep{\fill}}ld{2.0}d{1.1}cd{3.1}d{3.1}d{3.1}d{2.1}d{2.1}d{2.1}@{}}
\hline
&&& \multicolumn{7}{c@{}}{$\bolds{W_7}$}\\ [-6pt]
&&& \multicolumn{7}{c@{}}{\hrulefill}  \\
\multicolumn{1}{@{}l}{$\bolds{(n,p)}$} & \multicolumn{1}{c}{$\bolds{k}$} & \multicolumn{1}{c}{$\bolds{W_2}$} &
\multicolumn{1}{c}{$\bolds{\delta= }$} & \multicolumn{1}{c}{$\bolds{0.3}$} & \multicolumn{1}{c}{$\bolds{0.4}$} & \multicolumn{1}{c}{$\bolds{0.5}$} & \multicolumn{1}{c}{$\bolds{0.6}$} &
\multicolumn{1}{c}{$\bolds{0.7}$} & \multicolumn{1}{c}{$\bolds{0.8}$} \\
 \hline
$(60,40)$& 2 & 12.6 && 99.2 & 96.9 & 88.5 & 71.2 & 51.1 & 34.9  \\
& 4 & 15.9 && 99.4 & 98.4 & 92.9 & 75.4 & 54.1 & 39.3  \\
& 6 & 15.4 && 99.7 & 97.9 & 91.8 & 73.7 & 53 & 36.2  \\
& 8 & 16 && 99.6 & 98.2 & 90.9 & 70 & 47.5 & 33.3 \\
& 10 & 12.8 && 99.5 & 97.5 & 88.7 & 66.7 & 42.3 & 29.2  \\ [3pt]
$(120,80)$& 2 & 13.4 && 100 & 100 & 99.9 & 96.5 & 76.6 & 50.9  \\
& 4 & 13.2 && 100 & 100 & 100 & 97.5 & 79.2 & 54.3  \\
& 6 & 14.6 && 100 & 100 & 100 & 96.6 & 77.7 & 52.5  \\
& 8 & 15.2 && 100 & 100 & 100 & 96.5 & 76.2 & 49.3  \\
& 10 & 18.8 && 100 & 100 & 100 & 97.2 & 76.4 & 51.4  \\ [3pt]
$(60,10)$& $2$ & 91.6 && 100 & 100 & 99.9 & 99.7 & 99.8 & 98.7 \\
& $4$ & 84.9 && 100 & 100 & 99.9 & 99.5 & 98.4 & 96.9  \\
& $6$ & 83.4 && 100 & 100 & 100 & 99.4 & 98 & 96.4  \\
& $8$ & 80.8 && 100 & 100 & 99.9 & 99.3 & 97 & 95  \\
& $10$ & 83.1 && 100 & 100 & 99.9 & 98.7 & 97 & 94  \\[3pt]
$(120,160)$ & $2$ & 5.2 && 100 & 100 & 99.8 & 87.5 & 51.7 & 26.1  \\
& $4$ & 7.1 && 100 & 100 & 99.7 & 92 & 61.6 & 33.7  \\
& $6$ & 7.8 && 100 & 100 & 99.7 & 91.8 & 61.3 & 34.7  \\
& $8$ & 9.2 && 100 & 100 & 99.7 & 91.7 & 62.6 & 35.3  \\
& $10$ & 8.9 && 100 & 100 & 99.9 & 91.2 & 59.1 & 32.8  \\
\hline
\end{tabular*}
\end{table}
%

%s3 #&#
\section{Preliminaries and tools from Anderson and Zeitouni \texorpdfstring{\cite{AZ2009}}{[1]}}\label{sec3}

In this section, we will introduce some basic notions concerning
cumulants and some necessary results from \cite{AZ2009}. For some
positive integer $N$ and random variables $\xi_1,\ldots,\xi_N$, the
{\emph{joint cumulant}} $\mathbf{C}(\xi_1,\ldots,\xi_N)$ is
defined by
%
%e3.1 #&#
\begin{equation}
\mathbf{C}(\xi_1,\ldots,\xi_N)=\mathrm{i}^{-N}
\frac{\partial
^N}{\partial x_1\cdots\partial x_N}\log\mathbb{E}\exp\Biggl(\sum_{j=1}^N
\mathrm{i}x_j\xi_j\Biggr)\Bigg|_{x_1=\cdots=x_N=0}. \label{1115.1}
\end{equation}
It is straightforward to check via above definition that the following
properties hold:
\begin{longlist}[{P1}:]
\item[{P1}:] $\mathbf{C}(\xi_1,\ldots,\xi_N)$ \textit{is a symmetric
function of} $\xi_1,\ldots,\xi_N$.
\item[{P2}:] $\mathbf{C}(\xi_1,\ldots,\xi_N)$ \textit{is a multilinear
function of} $\xi_1,\ldots,\xi_N$.
\item[{P3}:] $\mathbf{C}(\xi_1,\ldots,\xi_N)=0$ \textit{if the variables}
$\xi_i,i\in[N]$ \textit{can be split into two groups} $\{\xi_i\}_{i\in
\mathbb
{S}_1}$ \textit{and} $\{\xi_i\}_{i\in\mathbb{S}_2}$ with $\mathbb{S}_1\cap
\mathbb{S}_2=\varnothing$ \textit{and} $\mathbb{S}_1\cup\mathbb{S}_2=[N]$ \textit{such
that the sigma field} $\sigma\{\xi_i\}_{i\in\mathbb{S}_1}$ \textit{is
independent of} $\sigma\{\xi_i\}_{i\in\mathbb{S}_2}$.
\end{longlist}

It is well known that the joint cumulant can be expressed in terms of
moments. To state this expression, we need some notions about set partition.
For some positive number $N$, let $L_N$ be the lattice consisting of
all the partitions of $[N]$. We say $\pi=\{B_1,\ldots, B_m\}\in L_N$ is
a {\emph{partition}} of the set $[N]$ if
\[
\varnothing\neq B_i\subset[N], i=1,\ldots,m,\qquad \bigcup
_{i=1}^mB_i=[N],\qquad B_i\cap
B_j= \varnothing, \mbox{ if } i\neq j.
\]
We say $B_i$'s are {\emph{blocks}} of $\pi$ and $m$ is the {\emph
{cardinality}} of $\pi$. We will also conventionally use the notation
$\#\pi$ to denote the cardinality of a partition $\pi$ all the way.
$L_N$ is a {\emph{poset}} (partially ordered set) ordered by set
inclusion. Specifically, given two partitions $\pi$ and $\sigma$ in
$L_N$, we say $\pi\leq\sigma$ (or $\pi$ is a {\emph{refinement}} of
$\sigma$) if every block of $\pi$ is contained in a block of $\sigma$.
Now given two partitions $\sigma_1,\sigma_2\in L_N$, with the above
order ``$\leq$'' we define $\sigma_1\vee\sigma_2$ to be the least
upper bound of $\sigma_1$ and~$\sigma_2$. For example, let $N=8$ and
\[
\sigma_1=\bigl\{\{1,2\},\{3,4,5\}, \{6\},\{7,8\}\bigr\},\qquad
\sigma_2=\bigl\{\{1,3\} ,\{2,5\},\{4\}, \{6,8\},\{7\}\bigr\}.
\]
Then we have $\sigma_1\vee\sigma_2=\{\{1,2,3,4,5\},\{6,7,8\}\}$.
With these notation, we have the following basic expression of joint
cumulant in terms of moments,
%
%e3.2 #&#
\begin{equation}
\mathbf{C}(\xi_1,\ldots,\xi_N)=\sum
_{\pi\in L_N}(-1)^{\#\pi-1}(\# \pi -1)!\mathbb{E}_\pi(
\xi_1,\ldots,\xi_N), \label{106.2}
\end{equation}
where
$
\mathbb{E}_\pi(\xi_1,\ldots,\xi_N)=\prod_{A\in\pi}\mathbb
{E}\prod_{i\in A}\xi_i$.\setcounter{footnote}{3}\footnote{Here, we remind that the partition $\pi$ in the notation
$\mathbb{E}_\pi(\cdot)$ takes effect on the positions of the components
of the vector, so does the notation $\mathbf{C}_\pi(\cdot)$ in (\ref
{925.1}). The reader should not confuse the positions with the
subscripts of the components of the random vector. For example, for
$\pi
=\{\{1,2\},\{3,4\}\}$, we have $\mathbb{E}_\pi
(Z_{j_2},Z_{j_3},Z_{j_4},Z_{j_1})=\mathbb{E}(Z_{j_2}Z_{j_3})\mathbb
{E}(Z_{j_4}Z_{j_1})$ rather than $\mathbb{E}(Z_{j_1}Z_{j_2})\mathbb
{E}(Z_{j_3}Z_{j_4})$.}
Especially, one has
\[
\mathbf{C}(\xi)=\mathbb{E}\xi, \qquad\mathbf{C}(\xi_1,\xi
_2)=\operatorname {Cov}(\xi_1,\xi_2).
\]

With the above concepts, we can now introduce the formula of joint
cumulants of $\operatorname{tr} S_n^{k_l},l=1,\ldots,r$ with $\sum_{l=1}^r
k_l=k$, derived in \cite{AZ2009}.
To this end, we need to specify two partitions $\pi_0,\pi_1\in
L_{2k}$ as
\begin{eqnarray*}
\pi_0&=&\bigl\{\{1,2\},\{3,4\},\ldots,\{2k-1,2k\}\bigr\},
\\
\pi_1&=&\bigl\{\{2,3\},\ldots, \{\mathbf{k}_1,1\},\{
\mathbf {k}_1+2,\mathbf {k}_1+3\},\ldots,
\\
&& \{\mathbf{k}_2,\mathbf{k}_1+1\},\ldots,\{\mathbf
{k}_{r-1}+2,\mathbf{k}_{r-1}+3\},\ldots,\{
\mathbf{k}_r,\mathbf {k}_{r-1}+1\}\bigr\},
\end{eqnarray*}
where $\mathbf{k}_i=2\sum_{j=1}^i k_j,i=1,\ldots,r$ and $k_j$'s are
nonnegative integers. Observe that $\#\pi_0\vee\pi_1=r$.
Now let $L_{N}^{2+}\subset L_N$ be the set consisting of those
partitions in which each block has cardinality larger than $2$. If each
block of a partition has cardinality~$2$, we call such a partition a
{\emph{perfect matching}}. For even $N$, we let $L_{N}^{2}$  ($\subset
L_N^{2+}$) be the set consisting of all perfect matchings. In the
sequel, we will also use the notation $L_{2k}^4$ to denote the set
consisting of the partitions of $[2k]$ containing one $4$-element block
and $(k-2)$ $2$-element blocks. Let $j_\alpha\in[n]$ for $\alpha\in
[N]$. We use the terminology in \cite{AZ2009} to call the index vector
$\mathbf{j}:=(j_1,\ldots, j_N)$ a $(n,N)$-{\emph{word}}. Moreover, we
say an $(n,N)$-word $\mathbf{j}$ is $\pi$-measurable for some partition
$\pi\in L_N$ when $j_\alpha=j_\beta$ if $\alpha,\beta$ are in the same
block of $\pi$.
Then by the discussions in \cite{AZ2009} (see Proposition~5.2 therein),
one has
%
%e3.3 #&#
\begin{eqnarray}
\label{925.1} &&\mathbf{C} \bigl(\operatorname{tr} S^{k_1}_n,
\ldots, \operatorname{tr} S^{k_r}_n\bigr)
\nonumber
\\[-8pt]
\\[-8pt]
\nonumber
&&\qquad = \mathop{\sum_{\pi\in L_{2k}^{2+}}}_{ \mathrm{s.t.}\ \# \pi
_0\vee\pi_1\vee\pi=1}
p^{-k+\# \pi_0\vee\pi} \mathop{\sum_{\mathbf{j}:(n,2k)\mbox{-}\mathrm{word}}}_{ \mathrm{s.t.}\ \mathbf{j}\ \mathrm{is}\ \pi_1\ \mathrm{measurable}}
\mathbf{C}_\pi(\mathbf{j}),
\end{eqnarray}
where
\[
\mathbf{C}_\pi(\mathbf{j}):=\mathbf{C}_\pi(Z_{j_1},
\ldots ,Z_{j_{2k}})=\prod_{A\in\pi} \mathbf{C}
\{Z_{j_\alpha}\}_{\alpha
\in A}.
\]
Equation~(\ref{925.1}) was derived by using M\"{o}bius inversion formula in
\cite
{AZ2009}. We remind that, here we switch the roles of the parameters
$p$ and $n$ in the setting of \cite{AZ2009}. Moreover, $B(\mathbf{j})$
therein is always $1$ for all $(n,2k)$-words $\mathbf{j}$ in our case.

Our aim is
to show that for $k_1,\ldots,k_r$ with $\sum_{j=1}^rk_j=k$,
%
%e3.4 #&#
\begin{eqnarray}
\label{1011.1} \mathbf{C} \bigl(\operatorname{tr} S^{k_1}_n,
\ldots, \operatorname{tr} S^{k_r}_n\bigr)= \cases{
o(1), &\quad $\mbox{if }r\geq3,$
\vspace*{2pt}\cr
\mbox{r.h.s. of }(\ref{1111.2})+o(1), &\quad$\mbox{if } r=2,$
\vspace*{2pt}\cr
\mbox{r.h.s. of }(\ref{1111.3}), &\quad $\mbox{if } r=1.$}
\end{eqnarray}
It is well known that (\ref{1011.1}) can imply Theorem~\ref{thmm.1111.1}
directly.

Apparently, by (\ref{925.1}) we see that, to prove (\ref{1011.1}), the
main task is to estimate the summation
%
%e3.5 #&#
\begin{equation}
\mathop{\sum_{\mathbf{j}: (n,2k)\mbox{-}\mathrm{word}}}_{
\mathrm{s.t.}\ \mathbf{j}\
\mathrm{is}\ \pi_1\ \mathrm{measurable}}
\mathbf{C}_\pi(\mathbf{j}) \label{1012.3}
\end{equation}
for various $\pi$. In order to deal with the estimation in the
counterpart of \cite{AZ2009}, the main assumption imposed on the
stationary sequence $\{Y_k\}_{k=-\infty}^{\infty}$ considered therein
is the so-called {\emph{joint cumulant summability}}
%
%e3.6 #&#
\begin{equation}
\sum_{j_1}\cdots\sum_{j_r}\bigl|
\mathbf{C}(Y_0,Y_{j_1},\ldots, Y_{j_r})\bigr|=O(1)\qquad
\mbox{for all } r\geq1. \label{1111.4}
\end{equation}
Actually, once joint cumulant summability held for random sequence $\{
Z_i\}_{i=1}^n$, one could obtain that the summation (\ref{1012.3}) was
bounded by $O(n^{\#\pi_1\vee\pi})$ in magnitude. However, the
verification of the joint cumulant summability for a general random
sequence is quite nontrivial. Without joint cumulant summability, we
will provide a weaker bound on (\ref{1012.3}) below directly; see
Proposition~\ref{lem.106.15} and Corollary~\ref{cor.107.20}. Such a
weaker bound will still make our proof strategy amenable. The following
proposition proved by Anderson and Zeitouni in \cite{AZ2009} will be
crucial to ensure that our weaker bound on (\ref{1012.3}) still works
well in the proof of (\ref{1011.1}).

%pr3.1 #&#
\begin{pro}[(Proposition~3.1, \cite{AZ2009})] \label{pro.1122.1}Let
$\Pi
\in L^{2+}_{2k}$ and $\Pi_0,\Pi_1\in L_{2k}^{2}$ for some positive
integer $k$. If $\#\Pi\vee\Pi_0\vee\Pi_1=1$ and $\#\Pi_0\vee\Pi_1=r$
for some positive integer $r\leq k$, one has
\[
\#\Pi_0\vee\Pi+\#\Pi_1\vee\Pi\leq(\#\Pi+1)
\mathbf{1}_{\{r=1\}
}+\min\{ \#\Pi+1, k+1-r/2\}\mathbf{1}_{\{r\geq2\}}.
\footnotemark[5]
\]
\end{pro}

\footnotetext[5]{One might note that in Proposition~3.1 of \cite{AZ2009},
the authors stated a slightly weaker bound $k+1-\lfloor r/2\rfloor$ in
formula (14) therein. However, according to the proof of Proposition~3.1 in \cite{AZ2009}, it is not difficult to see that one can improve
it to $k+1-r/2$. We appreciate Professor Greg W. Anderson's
confirmation on this.}

At the end of this section, we state some elementary properties of the
vector $\mathbf{Z}$ which will be used in the subsequent sections. We
summarize them as the following lemma whose proof will be stated in the
supplementary material \cite{BLPZ}.

%le3.2 #&#
\begin{lem} \label{lem.1011.2}Let $\mathbf{Z}=(Z_1,\ldots,Z_n)$ be the
random vector defined above. Let $m$ and $\alpha_i,i=1,\ldots,m$ be
fixed positive integers. We have the following properties of $\mathbf
{Z}$:
\begin{longlist}[(ii)]
\item[(i)] If $\sum_{i=1}^m\alpha_i$ is odd,
%
%e3.7 #&#
\begin{equation}
\mathbb{E}\bigl(Z_1^{\alpha_1}\cdots Z_m^{\alpha_m}
\bigr)=0. \label{926.2}
\end{equation}

\item[(ii)] If $\sum_{i=1}^m\alpha_i$ is even,
%
%e3.8 #&#
\begin{equation}
\mathbb{E}\bigl(Z_1^{\alpha_1}\cdots Z_m^{\alpha_m}
\bigr)=O\bigl(n^{-{\mathfrak {n}_o(\bolds{\alpha})}/{2}}\bigr), \label{1011.9}
\end{equation}
where $\bolds{\alpha}=(\alpha_1,\ldots,\alpha_m)$ and
$
\mathfrak{n}_o(\bolds{\alpha})=\#\{i\in[m]:\alpha_i\mbox{ is
odd}\}$.
Moreover, by symmetry, \textup{(i)} and \textup{(ii)} still hold if we replace
$(Z_1,\ldots,Z_m)$ by $(Z_{l_1},\ldots,Z_{l_m})$ with any mutually
distinct indices $l_1,\ldots, l_m\in[n]$.
\end{longlist}
\end{lem}

Note that (\ref{106.2}) together with (i) of Lemma~\ref{lem.1011.2}
implies that
%
%e3.9 #&#
\begin{equation}
\mathbf{C}(Z_{j_1},\ldots, Z_{j_r})=0\qquad \mbox{if } r\mbox{ is
odd}.
\end{equation}
Consequently, it suffices to consider those partition $\pi$ in which
each block has even cardinality. We denote $L_{2k}^{\operatorname
{even}}\ (\subset L_{2k}^{2+}$) to be the set of such partitions. Then we
can rewrite (\ref{925.1}) as
%
%e3.10 #&#
\begin{eqnarray}
\label{1012.1} &&\mathbf{C} \bigl(\operatorname{tr} S^{k_1}_n,
\ldots, \operatorname{tr} S^{k_r}_n\bigr)
\nonumber
\\[-8pt]
\\[-8pt]
\nonumber
&&\qquad = \mathop{\sum_{\pi\in L_{2k}^{\operatorname{even}}}}_{ \mathrm{s.t.}\ \#
\pi_0\vee\pi_1\vee\pi=1}
p^{-k+\# \pi_0\vee\pi} \mathop{\sum_{\mathbf{j}:
(n,2k)\mbox{-}\mathrm{word}}}_{ \mathrm{s.t.}\ \mathbf{j}\ \mathrm{is}\ \pi_1\ \mathrm{measurable}}
\mathbf{C}_\pi(\mathbf{j}).
\end{eqnarray}

%s4 #&#
\section{A rough bound on the summation (\texorpdfstring{\protect\ref{1012.3}}{3.5})}\label{sec4}
%s4.1 #&#
\subsection{Factorization of the summation}\label{sec4.1} As mentioned in the last
section, instead of checking the property of joint cumulant summability,\setcounter{footnote}{5}\footnote{As
mentioned above, once joint cumulant summability held, the
magnitude of (\ref{1012.3}) could be bounded by $O(n^{\#\pi_1\vee\pi
})$. To see this, one can refer to Proposition~6.1 of \cite{AZ2009} by
setting $b=p$ therein and switching the role of $p$ by $n$ to adapt to
our notation.} we will try to provide a rough bound on the summation
(\ref{1012.3})
for any given $\pi\in L_{2k}^{\operatorname{even}}$ directly. Now we
observe by
definition that
\[
\mathbf{C}_\pi(\mathbf{j})=\prod_{A\in\pi}
\mathbf{C}\{ Z_{j_\alpha}\} _{\alpha\in A}=\prod
_{B\in\pi_1\vee\pi} \mathop{\prod_{A\in\pi}}_{
\mathrm{ s.t.}\ A\subset B}
\mathbf{C}\{Z_{j_\alpha}\}_{\alpha\in A}.
\]
For simplicity, we denote $q:=\# \pi_1\vee\pi$ and index the blocks in
$\#\pi_1\vee\pi$ in any fixed order as $B_i, i\in[q]$ for convenience.
We set
\[
\mathfrak{b}_i:=\frac{\# B_i}{2},\qquad m_i:=\#\{A\in\pi: A
\subset B_i\}, \qquad i\in[q],
\]
and index the subsets of $B_i$ in $\pi$ in any fixed order by
$A_i^{(\beta)},\beta\in[m_i]$. Moreover, we define
\[
\mathfrak{a}_i(\beta):=\frac{\#A_i^{(\beta)}}{2}.
\]
Then we can write
%
%e4.1 #&#
\begin{equation}
\mathbf{C}_\pi(\mathbf{j})=\prod_{i=1}^q
\prod_{\beta
=1}^{m_i}\mathbf {C}
\{Z_{j_\alpha}\}_{\alpha\in A_i^{(\beta)}}. \label{1012.2}
\end{equation}

To simplify the notation, we use $\mathcal{J}:=\mathcal{J}(\pi_1,n,k)$
to represent the set consisting of all $\pi_1$-measurable $(n,2k)$-word
$\mathbf{j}$. Fix $\alpha\in[2k]$, we can view $\alpha: \mathcal
{J}\rightarrow[n]$ defined by $\alpha(\mathbf{j}):=j_\alpha$ as a
functional on $\mathcal{J}$. Then it is apparent that if $j_\alpha
\equiv j_\beta$ [i.e., $\alpha(\mathbf{j})\equiv\beta(\mathbf
{j})$] on
$\mathcal{J}$ for fixed $\alpha,\beta\in[2k]$, one has the fact that
$\alpha$ and $\beta$ are in the same block of $\pi_1$. We have the
following lemma on the properties of the triple $(\mathcal{J},\pi,\pi
_1)$ whose proof will be put in the supplementary material \cite{BLPZ}.

%le4.1 #&#
\begin{lem} \label{lem.105.1} Regarding $\{j_\alpha\}_{\alpha\in[2k]}$
as a class of $2k$ functionals on $\mathcal{J}$, we have:
\begin{longlist}[(ii)]
\item[(i)] Given $i\in[q]$, for any $\alpha\in B_i$, there exists exactly one
$\gamma\in B_i$ such that $\alpha\neq\gamma$ but $j_\alpha\equiv
j_\gamma$ on $\mathcal{J}$.

\item[(ii)] Given $i\in[q]$ with $m_i\geq2$, for any proper subset $
P\subset
[m_i]$, there exists at least one $\alpha_1\in\bigcup_{\beta\in
P}A_i^{(\beta)}$, for any other $\alpha_2\in\bigcup_{\beta\in
P}A_i^{(\beta)}$ one has $j_{\alpha_1}\not\equiv j_{\alpha_2}$ on
$\mathcal{J}$.
\end{longlist}
\end{lem}

For convenience, we use the notation
%
%e4.2 #&#
\begin{equation}
\mathbf{j}|_B=(j_\alpha)_{\alpha\in B} \label{1029.1}
\end{equation}
for any $B\subset[2k]$. Analogously, for any partition $\sigma$ and
$D\subset[2k]$ we will use the notation
\[
\sigma|_D=\{E\in\sigma:E\subset D\},
\]
which can be viewed as the partition $\sigma$ restricted on the set
$D$. Actually, Lemma~\ref{lem.105.1} is just a direct consequence of
the fact that $\pi_1|{B_i}$ is a perfect matching and
\[
\#\pi|_{B_i}\vee\pi_1|_{B_i}=\#(\pi\vee
\pi_1)|_{B_i}=1, \qquad i\in[q].
\]
Analogous to the notation $\mathcal{J}(\pi_1, n,k)$, we denote
$\mathcal{J}(\pi_1|_{B_i}, n, \mathfrak{b}_i)$ to be the set consisting
of all $\pi_1|_{B_i}$-measurable words $\mathbf{j}|_{B_i}$.
By (\ref{1012.2}) and (i) of Lemma~\ref{lem.105.1}, we see that for
given $\pi$
%
%e4.3 #&#
\begin{equation}
(\ref{1012.3}):=\sum_{\mathbf{j}\in\mathcal{J}} \mathbf{C}_\pi
(\mathbf {j})=\prod_{i=1}^q \sum
_{\mathbf{j}|_{B_i}\in\mathcal{J}(\pi_1|_{B_i},
n, \mathfrak{b}_i)}\mathbf{C}_{\pi|_{B_i}}(
\mathbf{j}|_{B_i}) .\label{1024.3}
\end{equation}
Hence, to estimate (\ref{1012.3}), it suffices to provide a bound on
the quantity
\[
\sum_{\mathbf{j}|_{B_i}\in\mathcal{J}(\pi_1|_{B_i}, n, \mathfrak
{b}_i)}\mathbf{C}_{\pi|_{B_i}}(
\mathbf{j}|_{B_i})=\sum_{\mathbf
{j}|_{B_i}\in\mathcal{J}(\pi_1|_{B_i}, n, \mathfrak{b}_i)}\prod
_{\beta
=1}^{m_i}\mathbf{C}\{Z_{j_\alpha}
\}_{\alpha\in A_i^{(\beta)}}
\]
for all fixed $i\in[q]$.
For simplicity, we discard the subscript $i$ in the discussion below.
Thus, we will use the notation $B,A^{(\beta)}, \mathfrak{a}(\beta),
\mathfrak{b}, m$ to replace $B_i,A_i^{(\beta)},\mathfrak{a}_i(\beta),
\mathfrak{b}_i, m_i$ temporarily. Our main technical result is the
following crucial proposition.

%pr4.2 #&#
\begin{pro} \label{lem.106.15}With the above notation, we have
%
%e4.4 #&#
\begin{equation}
\sum_{\mathbf{j}|_{B}\in\mathcal{J}(\pi_1|_{B}, n, \mathfrak{b})}
 \prod_{\beta=1}^{m}\bigl|
\mathbf{C}\{Z_{j_\alpha}\}_{\alpha\in A^{(\beta
)}}\bigr|=O \bigl(n^{1+\sum_{\beta=1}^m(\mathfrak{a}(\beta)-2)\mathbf
{1}_{\{
\mathfrak{a}(\beta)\geq2\}}} \bigr).
\label{107.3}
\end{equation}
\end{pro}

%re4.3 #&#
\begin{rem} Here, we draw the attention that by (i) of Lemma~\ref
{lem.105.1}, the indices in $\mathbf{j}|_B$ are equivalent in pairs on
$\mathcal{J}(\pi_1|_B,n,\mathfrak{b})$. Under this constraint, the
number of free $\mathbf{j}$ indices in the summation on the left-hand
side of (\ref{107.3}) is actually $\mathfrak{b}=\sum_{i=1}^m
\mathfrak
{a}(\beta)$.
\end{rem}

Now we use $\mathfrak{n}_\gamma(\pi)$ to represent the number of blocks
in $\pi$ whose cardinalities are $2\gamma$. We have the following
corollary whose proof follows from Proposition~\ref{lem.106.15} and
(\ref{1024.3}) directly.

%co4.4 #&#
\begin{cor} \label{cor.107.20}With the above notation, we have
%
%e4.5 #&#
\begin{equation}
\sum_{\mathbf{j}\in\mathcal{J}}\bigl|\mathbf{C}_\pi(\mathbf{j})\bigr|=
O\bigl(n^{\#\pi
_1\vee\pi+\sum_{\gamma\geq2}(\gamma-2)\mathfrak{n}_\gamma(\pi
)}\bigr). \label{1015.4}
\end{equation}
\end{cor}

Our main task in this section is to prove Proposition~\ref{lem.106.15}.
The tedious proof will be given in the remaining part of this section
which will be further split into several subsections. In Sections~\ref
{sec4.2}--\ref{sec4.4}, we will provide some preliminary results for
our final
evaluation scheme. The formal proof of Proposition~\ref{lem.106.15}
will be stated in Section~\ref{sec4.5}.
%s4.2 #&#
\subsection{Bounds on joint cumulants}\label{sec4.2} In this subsection, we will
provide some bounds on single joint cumulants with variables from
$\mathbf{Z}$. Such bounds will help us to reduce all these joint
cumulants to some products of $2$-element cumulants which are more
friendly for the subsequent combinatorial enumeration. Let $s,t$ be
fixed nonnegative integers.
Now for $l_1,\ldots,l_s, h_1,\ldots, h_{2t}\in[n]$, we denote the vectors
\[
\mathbf{l}:=(l_1,l_1,l_2,l_2,
\ldots,l_s,l_s), \qquad \mathbf {h}:=(h_1,
\ldots,h_{2t})
\]
and we will use {\emph{$\mathbf{l}$ index}} (resp., {\emph{$\mathbf{h}$
index}}) to refer to $l_i,i\in[s]$ (resp., $h_i,i\in[2t]$). For
simplicity, we will employ the notation
\[
\widehat{\mathbf{l}\mathbf{h}}=(l_1,l_1,l_2,l_2,
\ldots,l_s,l_s, h_1,\ldots,h_{2t}),
\]
which is the concatenation of $\mathbf{l}$ and $\mathbf{h}$.
Note that the $\mathbf{l}$ indices appear in pairs. And in the sequel,
we will use the notation
%
%e4.6 #&#
\begin{equation}
\mathbf{C}(\widehat{\mathbf{l}\mathbf{h}}):=\mathbf {C}(Z_{l_1},Z_{l_1},
\ldots, Z_{l_s},Z_{l_s}, Z_{h_1},\ldots,
Z_{h_{2t}})\label{108.1}
\end{equation}
and
\[
\mathbf{C}(\mathbf{h}):=\mathbf{C}(Z_{h_1},\ldots, Z_{h_{2t}})
\]
to highlight the index sequence. In this manner, for any partition
$\tilde{\pi}\in L_{2s+2t}$, we denote
\[
\mathbb{E}_{\tilde{\pi}}(\widehat{\mathbf{l}\mathbf{h}}):=\mathbb
{E}_{\tilde{\pi}}(Z_{l_1},Z_{l_1},\ldots,
Z_{l_s},Z_{l_s}, Z_{h_1},\ldots, Z_{h_{2t}}).
\]
Analogously, for any partition $\tilde{\sigma}\in L_{2t}$, we set
\[
\mathbb{E}_{\tilde{\sigma}}(\mathbf{h}):=\mathbb{E}_{\tilde
{\sigma
}}(Z_{h_1},Z_{h_2},
\ldots, Z_{h_{2t}}).
\]
We remind here the aforementioned convention that the partition $\tilde
{\pi}$ in the notation $\mathbb{E}_{\tilde{\pi}}(\cdot)$ takes effect
on the positions of components of $\widehat{\mathbf{l}\mathbf{h}}$.

To prove Proposition~\ref{lem.106.15}, we will need the following two
lemmas. Lemma~\ref{lem.1117.10} gives us an explicit order on the
magnitude of $4$-element cumulant, whilst Lemma~\ref{lem.106.16}
provides some rough bounds on the cumulants with more than $4$ elements.

%le4.5 #&#
\begin{lem} \label{lem.1117.10} Suppose that $\mathbf
{h}=(h_1,h_2,h_3,h_4)$. Let $d(\mathbf{h})$ be the number of the
distinct values in $\{h_1,h_2,h_3,h_4\}$. We have
%
%e4.7 #&#
\begin{equation}
\mathbf{C}(\mathbf{h})=O\bigl(n^{-d(\mathbf{h})+1}\bigr). \label{106.20}
\end{equation}
\end{lem}

The proof of Lemma~\ref{lem.1117.10} will be stated in the
supplementary material \cite{BLPZ}. From Lemma~\ref{lem.1117.10}, we
can get the following consequences. We see that if there exists a
perfect matching $\sigma=\{A_1,A_2\}$ of $\{1,2,3,4\}$ such that $\{
h_i\}_{i\in A_1}\cap\{h_i\}_{i\in A_2}=\varnothing$, then by using Lemma~\ref{lem.1117.10} and (\ref{1113.20}) we can get
%
%e4.8 #&#
\begin{equation}
\bigl|\mathbf{C}(\mathbf{h})\bigr|\leq O\bigl(n^{-1}\bigr)\bigl|\mathbb{E}_{\sigma}(
\mathbf {h})\bigr|. \label{107.12}
\end{equation}
If there is no such perfect matching, we have
%
%e4.9 #&#
\begin{equation}
\bigl|\mathbf{C}(\mathbf{h})\bigr|\leq C\bigl|\mathbb{E}_{\sigma}(\mathbf
{h})\bigr|\label{107.13}
\end{equation}
for any perfect matching $\sigma\in L_4^2$ with some positive constant
$C$. Note that the second case occurs if and only if three or four of
the indices $h_1,h_2,h_3,h_4$ take the same value.

The following Lemma~\ref{lem.106.16} provides some crucial bounds on
the cumulants with more than $4$ underlying elements. The proof of this
lemma will also be stated in the supplementary material \cite{BLPZ}.

%le4.6 #&#
\begin{lem} \label{lem.106.16}Under the above notation, we have the
following bounds on the joint cumulant $\mathbf{C}(\widehat{\mathbf
{l}\mathbf{h}})$.
\begin{longlist}[(iii)]
\item[(i)] (Crude bound) When $s+t\geq3$ and $t\geq1$,
%
%e4.10 #&#
\begin{equation}
\bigl|\mathbf{C}(\widehat{\mathbf{l}\mathbf{h}})\bigr|\leq C\sum
_{\sigma\in
L_{2t}^2} \bigl|\mathbb{E}_\sigma(\mathbf{h})\bigr| \label{141241}
\end{equation}
for some positive constant $C$.

\item[(ii)] When $s\geq2$, $t=1$, $l_1,\ldots, l_s$ are mutually
distinct and distinct from $h_1,h_2$, we have
%
%e4.11 #&#
\begin{equation}
\bigl|\mathbf{C}(\widehat{\mathbf{l}\mathbf{h}})\bigr|\leq O\bigl(n^{-2}\bigr).
\label{141242}
\end{equation}
\item[(iii)] When $s\geq2$, $t=1$, $l_1=l_2$ and $l_2,\ldots, l_s$
are mutually distinct and distinct from $h_1,h_2$, we have
%
%e4.12 #&#
\begin{equation}
\bigl|\mathbf{C}(\widehat{\mathbf{l}\mathbf{h}})\bigr|=O\bigl(n^{-1}\bigr)\bigl|
\mathbb {E}(Z_{h_1}Z_{h_2})\bigr|. \label{141243}
\end{equation}
\end{longlist}
\end{lem}

At the end of this subsection, we need to clear up a potential
confusion which may occur when we use Lemma~\ref{lem.106.16} in the
sequel. Given an index sequence, for example, $(1,1,2,2,3,3,4,4)$, we
consider to use Lemma~\ref{lem.106.16} to bound the corresponding
cumulant $\mathbf{C}(Z_1,Z_1,\ldots,Z_4,Z_4)$. Obviously, we can adopt
(ii) of Lemma~\ref{lem.106.16} by setting $l_1=1$, $l_2=2$, $l_3=3$ and
$h_1=h_2=4$. Thus, $s=3$ and $t=1$. However, we can also say that
$s=t=2$ such that $l_1=1$, $l_2=2$ while $h_1=h_2=3$ and $h_3=h_4=4$.
We can even say that $s=0,t=4$ such that $h_{2i-1}=h_{2i}=i,i=1,\ldots
,4$. That means the determination of $\mathbf{l}$ and $\mathbf{h}$
indices as well as $s$ and $t$ are not substantially important.
Actually in any viewpoint listed above, we can employ (ii) of Lemma~\ref
{lem.106.16}. We state Lemma~\ref{lem.106.16} with $\mathbf{l}$ and
$\mathbf{h}$ in this way in order to simplify the presentation.
However, when we use Lemma~\ref{lem.106.16}, we only need to check
which case of (i)--(iii) is applicable to the given index sequence.
Moreover, we can also represent the bounds for (ii) and (iii) in the
form of the right-hand side of (\ref{141241}). In case (ii), obviously,
we can find a perfect matching $\tilde{\sigma}\in L_{2s+2}^2$ such that
\[
\mathbf{C}_{\tilde{\sigma}}(\widehat{\mathbf{l}\mathbf{h}})=\prod
_{i=1}^s\mathbf{C}(Z_{l_i},Z_{l_i})
\cdot\mathbf {C}(Z_{h_1},Z_{h_2})=\mathbb{E}(Z_{h_1},Z_{h_2}).
\]
By (\ref{1113.20}) and (\ref{141242}), we observe that in case (ii),
%
%e4.13 #&#
\begin{eqnarray}
\label{141244}\bigl |\mathbf{C}(\widehat{\mathbf{l}\mathbf{h}})\bigr|&\leq& O
\bigl(n^{-1}\bigr)\bigl|\mathbb {E}(Z_{h_1},Z_{h_2})\bigr|=O
\bigl(n^{-1}\bigr)\bigl|\mathbb{E}_{\tilde{\sigma}}(\widehat {\mathbf{l}
\mathbf{h}})\bigr|
\nonumber
\\[-8pt]
\\[-8pt]
\nonumber
&\leq& O\bigl(n^{-1}\bigr)\sum_{\sigma\in L^2_{2s+2}}\bigl|
\mathbb{E}_\sigma (\widehat {\mathbf{l}\mathbf{h}})\bigr|.
\end{eqnarray}
Note that the above bound is not as strong as (\ref{141242}) when
$h_1=h_2$. Analogously, it is easy to check that (\ref{141244}) also
holds in case (iii) of Lemma~\ref{lem.106.16}.
%s4.3 #&#
\subsection{Cyclic product of $2$-element cumulants}\label{sec4.3}
In this subsection, we introduce the concept of cyclic product of
$2$-element cumulant ({\emph{cycle}} in short) and the summation of
this kind of products over involved components of $\mathbf{j}$ words.
Such products will serve as canonical factors in the subsequent
discussion on the whole product $\mathbf{C}_{\pi|_B}(\mathbf{j}|_B)$ in
Proposition~\ref{lem.106.15}. Let $\ell$ be some positive integer and
$\sigma\in L_{2\ell}^2$. As above, we use the notation $\mathcal
{J}(\sigma,n,\ell)$ to denote the set of all $\sigma$-measurable
$(n,2\ell)$-words $\mathbf{j}$. Moreover, let $\sigma_0\in L_{2\ell
}^2$. Note that $\mathbf{C}_{\sigma_0}(\mathbf{j})$ is a product of
$\ell$ $2$-element cumulants. Now we define the concept of {\emph
{cycle}} (with respect to $\sigma$) as follows.

%de4.7 #&#
\begin{defi} [(Cycle)]\label{defi.1116.2} Let $\sigma_0, \sigma\in
L_{2\ell}^2$ and $\mathbf{j}\in\mathcal{J}(\sigma,n,\ell)$. We
call the
cumulant product $\mathbf{C}_{\sigma_0}(\mathbf{j})$ a cycle with
respect to $\sigma$ if $\#\sigma\vee\sigma_0=1$.
\end{defi}

%re4.8 #&#
\begin{rem} Note that actually for $\mathbf{j}\in\mathcal{J}(\sigma
,n,\ell)$, whether $\mathbf{C}_{\sigma_0}(\mathbf{j})$ is a cycle (with
respect to $\sigma$) only depends on the perfect matchings $\sigma$ and
$\sigma_0$ but not on the choice of $\mathbf{j}$. However, the
magnitude of a cycle $\mathbf{C}_{\sigma_0}(\mathbf{j})$ does depend on
the choice of the word $\mathbf{j}$. See Lemma~\ref{lem.1105.2} below.
\end{rem}

\begin{figure}[b]

\includegraphics{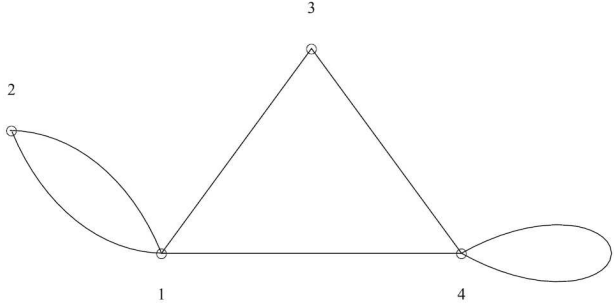}

\caption{Let $\sigma=\{\{1,2\},\{3,4\},\{5,6\},\{7,8\},\{9,10\},\{
11,12\}\}$ and $\sigma_0=\{\{2,3\},\{4,5\},\break \{6,7\}, \{8,9\}, \{10,11\},\{
12,1\}\}$. We take the $\sigma$-measurable word $\mathbf{j}$ to be
$(1,1,2,2,1,1,3,3,\break 4,4,4,4)$. Then the corresponding graph $G(\mathbf
{j})$ for the cycle $\mathbf{C}_{\sigma_0}(\mathbf{j})$ is as above.}
\label{fig1}
\end{figure}

We can illustrate the definition of {\emph{cycle}} in the following
more detailed way through which the meaning of such a nomenclature can
be evoked. Provided that $\#\sigma\vee\sigma_0=1$, it is not difficult
to see that there exists a permutation $\varepsilon$ of $[2\ell]$
such that
\[
\sigma_0=\bigl\{\bigl\{\varepsilon(2\alpha-1),\varepsilon(2\alpha)
\bigr\}\bigr\} _{\alpha
=1}^\ell,\qquad \sigma=\bigl\{\bigl\{\varepsilon(2
\alpha),\varepsilon(2\alpha+1)\bigr\} \bigr\}_{\alpha=1}^\ell
\]
in which we made the convention of $\varepsilon(2\ell+1)=\varepsilon
(1)$. Then for all $\mathbf{j}\in\mathcal{J}(\sigma,n,\ell)$,
$\mathbf
{C}_{\sigma_0}(\mathbf{j})$ can be written as a product of 2-element
cumulants whose indices form a cycle in the sense that if we regard
$V(\mathbf{j}):=\bigcup_{\alpha=1}^{2\ell}\{j_\alpha\}$ as the set of
vertices and $E(\mathbf{j}):=\bigsqcup_{\alpha=1}^{\ell}\{\overline
{j_{\varepsilon(2\alpha-1)}j_{\varepsilon(2\alpha)}}\}$ as the set of
edges then the multigraph $G(\mathbf{j})=(V(\mathbf{j}),E(\mathbf{j}))$
is a cycle (i.e., closed walk). Here, the notation $\cup$ is the common
union while $\sqcup$ is the disjoint union. The reader is recommended
to take a look at Figure~\ref{fig1} for understanding the definition of
a cycle. In this manner, we will also say that the word $\mathbf
{j}=(j_1,\ldots,j_{2\ell})$ forms a cycle under $\sigma_0$ with respect
to $\sigma$.

%f1 #&#

%de4.9 #&#
\begin{defi}[(Product of $\mathfrak{m}$ cycles)] Given $\sigma,\sigma
_0\in L_{2\ell}^2$, if $\#\sigma_0\vee\sigma=\mathfrak{m}$ for some
positive integer $\mathfrak{m}\geq2$, we say that $\mathbf
{C}_{\sigma
_0}(\mathbf{j})$ with $\mathbf{j}\in\mathcal{J}(\sigma,n,\ell)$
is a
product of $\mathfrak{m}$ cycles with respect to $\sigma$. Actually, if
$\sigma_0\vee\sigma=\{D_1,\ldots,D_\mathfrak{m}\}$, then obviously
$\mathbf{C}_{\sigma_0|_{D_i}}(\mathbf{j}|_{D_i})$ with $\mathbf
{j}\in
\mathcal{J}(\sigma,n,\ell)$ is a cycle with respect to $\sigma|_{D_i}$.
\end{defi}

%re4.10 #&#
\begin{rem}Note that the definition of {\emph{product of $\mathfrak{m}$
cycles}} also only depends on $\sigma_0$ and $\sigma$. In some sense,
unlike the concept of {\emph{cycle}}, the {\emph{product of
$\mathfrak
{m}$ cycles}} is more like a common phrase rather than a new
terminology. However, we still raise it as an independent concept to
emphasize the relationship between $\sigma_0$ and $\sigma$.
\end{rem}

\begin{figure}[b]

\includegraphics{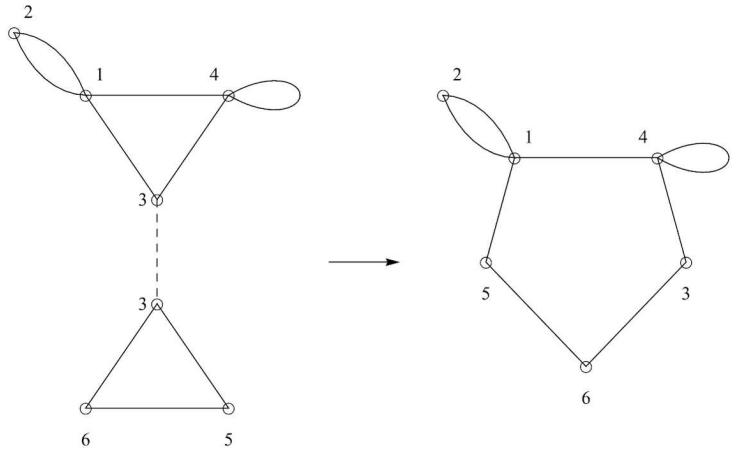}

\caption{Assume that $\sigma=\{\{2i-1,2i\}:i=1,\ldots,9\}$ and
$\sigma
_0=\{\{2,3\},\{4,5\},\{6,7\},\break  \{8,9\},\{10,11\},\{12,1\},\{13,18\},\{
14,15\},\{16,17\}\}$. We take the word $\mathbf{j}$ to be
$(1,1,2,2,\break 1,1,3,3,4,4,4,4,3,3,5,5,6,6)$. Then the figure on the
left-hand side above is $G(\mathbf{j})$ for the product of 2 cycles
$\mathbf{C}_{\sigma_0}(\mathbf{j})$. Now we fix a way to merge
$\mathbf
{C}_{\sigma_0}(\mathbf{j})$ by taking $\tilde{\sigma}_0=\{\{2,3\},\{
4,5\},\{6,15\},\{8,9\},\{10,11\},\{12,1\},\{13,18\},\{16,17\}\}$, then
the figure on the right-hand side above is corresponding to the merged
cycle.} \label{fig2}
\end{figure}

Actually, for some specific $\mathbf{j}$, the graphical illustration as
that in the single cycle case may not evoke the name of {\emph{product
of $\mathfrak{m}$ cycles}} any more since different graphs
corresponding to different cycles may have coincident vertices, and
thus these cycles will be tied together if we define $G(\mathbf
{j})=(V(\mathbf{j}),E(\mathbf{j}))$ as above. To avoid this confusion,
we can draw $\mathfrak{m}$ cycles separately and view the disjoint
union of these components as the graph corresponding to the product of
$\mathfrak{m}$ cycles. Moreover, we can use the dash line to connect
coincident indices in different components. One can see the left half
of Figure~\ref{fig2}, for example.
Actually, when there is some coincidence between indices from different
components, we will introduce a {\emph{merge operation}} to reduce the
number of cycles in the product later. Before commencing this issue, we
will raise a fact on the summation of single cycles.
Now we have the following lemma whose proof will be stated in the
supplementary material \cite{BLPZ}.

%f2 #&#

%le4.11 #&#
\begin{lem} \label{lem.1105.2}Let $\ell$ be a fixed positive integer,
and $\sigma_0,\sigma\in L_{2\ell}^2$ such that $\#\sigma_0\vee
\sigma
=1$. Assume that $\mathbf{j}\in\mathcal{J}(\sigma,n,\ell)$, and thus
$\mathbf{C}_{\sigma_0}(\mathbf{j})$ is a cycle with respect to~$\sigma
$. We have
%
%e4.14 #&#
\begin{equation}
\bigl|\mathbf{C}_{\sigma_0}(\mathbf{j})\bigr|= O\bigl(n^{-d(\mathbf{j})\mathbf
{1}_{\{
d(\mathbf{j})\geq2\}}}\bigr),
\label{1105.1}
\end{equation}
where $d(\mathbf{j})$ represents the number of distinct values in the
collection $\{j_\alpha\}_{\alpha=1}^{2\ell}$.
\end{lem}

In the sequel, we call a cycle containing at least one factor $\mathbf
{C}(Z_{j_{\varepsilon(2\alpha-1)}},\break Z_{j_{\varepsilon(2\alpha)}})$ with
$j_{\varepsilon(2\alpha-1)}\neq j_{\varepsilon(2\alpha)} $ as
{\emph
{in-homogeneous cycle}}. Otherwise, we call it {\emph{homogeneous
cycle}}. With the above graphical language, a homogeneous cycle only
has a single vertex and all its edges are self-loops. By contrast, an
in-homogeneous cycle has at least two vertices. Following from Lemma~\ref{lem.1105.2}, we have
the following.

%co4.12 #&#
\begin{cor} \label{cor.1105.4} For any given positive integer $\ell$,
and $\sigma_0,\sigma\in L_{2\ell}^2$ such that $\#\sigma\vee\sigma
_0=1$, we have the following corollary:
%
%e4.15 #&#
\begin{equation}
\sum_{\mathbf{j}\in\mathcal{J}(\sigma,n,\ell)}\bigl|\mathbf {C}_{\sigma
_0}(
\mathbf{j})\bigr|=n+O(1). \label{1105.3}
\end{equation}
\end{cor}

\begin{pf}
By using Lemma~\ref{lem.1105.2}, the leading term of the left-hand side
of (\ref{1105.3}) comes from the homogeneous cycles. Obviously, the
total choice of homogeneous cycle is $n$. Moreover, we can see that the
total contribution of the in-homogeneous cycles is $O(1)$ by (\ref
{1105.1}). Thus, we obtain (\ref{1105.3}).
\end{pf}

Now we use Corollary~\ref{cor.1105.4} to prove a simple case of
Proposition~\ref{lem.106.15}. That is $\mathfrak{a}(\beta)=1$ for all
$\beta=1,\ldots,m$. Note that in this case, by Lemma~\ref{lem.105.1},
it is not difficult to see that $\prod_{\beta=1}^{m}\mathbf{C}\{
Z_{j_\alpha}\}_{\alpha\in A^{(\beta)}}$ is a cycle for $\mathbf
{j}|_{B}\in\mathcal{J}(\pi_1|_{B}, n, \mathfrak{b})$ since $\#\pi
|_B\vee\pi_1|_B=1$. Hence, one has
%
%e4.16 #&#
\begin{equation}
\sum_{\mathbf{j}|_{B}\in\mathcal{J}(\pi_1|_{B}, n, \mathfrak
{b})}\prod_{\beta=1}^{m}\bigl|
\mathbf{C}\{Z_{j_\alpha}\}_{\alpha\in A^{(\beta
)}}\bigr|=n+O(1). \label{1129.5}
\end{equation}

We conclude this subsection by introducing the concept of {\emph{merge
operation}} toward the product of $\mathfrak{m}$ cycles when at least
two cycles in this product have some coincident indices. Now note that
if $\#\sigma\vee\sigma_0=2$, we see that for $\mathbf{j}\in
\mathcal
{J}(\sigma,n,\ell)$, $\mathbf{C}_{\sigma_0}(\mathbf{j})$ is a product
of two cycles by definition. In other words, there exist some
permutation $\varepsilon$ of $[2\ell]$ and $\ell_1\in[\ell]$ such that
\[
\sigma_0=\bigl\{\bigl\{\varepsilon(2\alpha-1),\varepsilon(2\alpha)
\bigr\}\bigr\} _{\alpha
=1}^\ell
\]
and
\begin{eqnarray*}
\sigma&=&\bigl\{\bigl\{\varepsilon(2),\varepsilon(3)\bigr\},\ldots,\bigl\{
\varepsilon (2\ell _1),\varepsilon(1)\bigr\},\bigl\{\varepsilon(2
\ell_1+2),\varepsilon(2\ell _1+3)\bigr\} ,\ldots,\\
&&{} \bigl
\{\varepsilon(2\ell),\varepsilon(2\ell_1+1)\bigr\}\bigr\}.
\end{eqnarray*}
Therefore, we have
\[
\mathbf{C}_{\sigma_0}(\mathbf{j})=\prod_{\alpha=1}^{\ell
_1}
\mathbf {C}(Z_{j_{\varepsilon(2\alpha-1)}}, Z_{j_{\varepsilon(2\alpha
)}})\prod
_{\alpha=\ell_1+1}^{\ell}\mathbf{C}(Z_{j_{\varepsilon(2\alpha-1)}},
Z_{j_{\varepsilon(2\alpha)}}).
\]
Now if for some specified $\mathbf{j}\in\mathcal{J}(\sigma,n,\ell)$,
there exist some $\beta\in[2\ell_1]$, $\gamma\in[2\ell]\setminus
[2\ell
_1]$ such that $j_{\varepsilon(\beta)}$ and $j_{\varepsilon(\gamma)}$
take the same value, we define the following {\emph{merge operation}}
for the two-cycle product $\mathbf{C}_{\sigma_0}(\mathbf{j})$. Without
loss of generality, we let $\beta=2\ell_1,\gamma=2\ell_1+1$. Then in
this case we have $j_{\varepsilon(1)}=j_{\varepsilon(2\ell
_1)}=j_{\varepsilon(2\ell_1+1)}=j_{\varepsilon(2\ell)}$ since
$\mathbf
{j}\in\mathcal{J}(\sigma,n,\ell)$. Now by (\ref{1113.20}) we see that
%
%e4.17 #&#
\begin{eqnarray}
\label{1116.1} &&\bigl|\mathbf{C}(Z_{j_{\varepsilon(2\ell_1-1)}},Z_{j_{\varepsilon
(2\ell
_1)}})
\mathbf{C}(Z_{j_{\varepsilon(2\ell_1+1)}},Z_{j_{\varepsilon
(2\ell
_1+2)}})\bigr|
\nonumber
\\[-8pt]
\\[-8pt]
\nonumber
&&\qquad \leq\bigl|\mathbf{C}(Z_{j_{\varepsilon(2\ell
_1-1)}},Z_{j_{\varepsilon(2\ell_1+2)}})\bigr|
\end{eqnarray}
when $j_{\varepsilon(2\ell_1)}=j_{\varepsilon(2\ell_1+1)}$. We set
\begin{eqnarray*}
\tilde{\sigma}_0&=&\bigl(\sigma_0\setminus\bigl\{
\bigl\{ \varepsilon (2\ell_1-1),\varepsilon(2\ell_1)
\bigr\},\bigl\{\varepsilon(2\ell _1+1),\varepsilon (2
\ell_1+2)\bigr\}\bigr\}\bigr)
\\
&&{}\cup\bigl\{\bigl\{\varepsilon(2\ell_1-1),\varepsilon(2
\ell_1+2)\bigr\}\bigr\}
\end{eqnarray*}
and
\[
\tilde{\sigma}=\bigl(\sigma\setminus\bigl\{\bigl\{\varepsilon(1),\varepsilon (2
\ell _1)\bigr\},\bigl\{\varepsilon(2\ell_1+1),
\varepsilon(2\ell)\bigr\}\bigr\}\bigr)\cup\bigl\{\bigl\{ \varepsilon(1),
\varepsilon(2\ell)\bigr\}\bigr\}.
\]
Then we obtain that
%
%e4.18 #&#
\begin{equation}
\label{1117.4} \mathbf{C}_{\tilde{\sigma}_0}(Z_{j(\varepsilon(1))}\cdots
Z_{j(\varepsilon(2\ell_1-1))},Z_{j(\varepsilon(2\ell_1+2))}, \ldots, Z_{j(\varepsilon(2\ell))})
\end{equation}
forms a cycle with respect to $\tilde{\sigma}$. We call (\ref{1117.4})
the {\emph{merged cycle}} of $\mathbf{C}_{\sigma_0}(\mathbf{j})$ (see
Figure~\ref{fig2}, e.g.).
Then by (\ref{1116.1}), we have
\begin{eqnarray*}
&&\mathop{\sum_{\mathbf{j}\in\mathcal{J}(\sigma,n,\ell)}}_{
\mathrm{subject\ to\ }j_{\varepsilon(2\ell_1)}=j_{\varepsilon(2\ell
_1+1)}}\bigl|
\mathbf{C}_{\sigma_0}(\mathbf{j})\bigr|
\\
&&\qquad\leq\sum_{\mathbf{j}\in\mathcal{J}(\tilde{\sigma},n, \ell-1)} \bigl|\mathbf{C}_{\tilde{\sigma}_0}(Z_{j(\varepsilon(1))}
\cdots Z_{j(\varepsilon(2\ell_1-1))},Z_{j(\varepsilon(2\ell_1+2))}, \ldots, Z_{j(\varepsilon(2\ell))})\bigr|
\\
&&\qquad =n+O(1).
\end{eqnarray*}
Obviously, the way to do the merge operation may be not unique when
there are more than one common value of the vertices from two different
cycles. In this case, we can just choose one way to do the merge
operation since in the sequel we only care about whether two cycles can
be merged but do not care about how to merge them. Analogously, in this
manner, when $\#\sigma\vee\sigma_0\geq3$, we can start from two cycles
and use the merge operation to merge them into one cycle once there
exists at least two indices (one from each cycle) taking the same
value, and then we can proceed this merge operation until there is no
pair of cycles can be merged into one.

%s4.4 #&#
\subsection{Classification of the relationships between indices}\label{sec4.4}
Note that the inequalities in Lemma~\ref{lem.106.16} rely on the
relationships between the underlying indices in the joint cumulants. In
order to use Lemmas \ref{lem.1117.10} and \ref{lem.106.16} in the proof
of Proposition~\ref{lem.106.15}, we will introduce some notation and
approach to classify the relationships between the indices (components
of $\mathbf{j}$).

At first, we introduce the concepts of {\emph{paired indices}} and
{\emph{unpaired indices}} as follows. Note that for any block $\{
d_1,d_2\}\in\pi_1$, we have $j_{d_1}\equiv j_{d_2}$ on $\mathcal{J}$
by definition. Now for each block $A^{(\beta)}\in\pi|_B$, we find out
all $\pi_1$'s blocks which are totally contained in $A^{(\beta)}$. We
denote the number of such blocks by $s(\beta)$. Specifically, we find
out all blocks $D^{(\beta)}_i:=\{d^{(\beta)}_{i1},d^{(\beta)}_{i2}\}
\in
\pi_1,i=1,\ldots, s(\beta)$ such that $\bigcup_{i=1}^{s(\beta
)}D^{(\beta
)}_i\subset A^{(\beta)}$. We then set
\[
l^{(\beta)}_i:=j_{d^{(\beta)}_{i1}}\equiv j_{d^{(\beta)}_{i2}}\qquad
\mbox{on }\mathcal{J}, i=1,\ldots, s(\beta).
\]
We call $l^{(\beta)}_i,i=1,\ldots, s(\beta)$ {\emph{paired indices}}
{\emph{from}} $A^{(\beta)}$ informally. The remaining indices
$j_\alpha
$ with $\alpha\in A^{(\beta)}\setminus\bigcup_{i=1}^{s(\beta
)}D^{(\beta
)}_i$ will be ordered (in any fixed order) and denoted by $h^{(\beta
)}_i, i=1,\ldots, 2t(\beta)$ which will be called as {\emph{unpaired
indices}} {\emph{from}} $A^{(\beta)}$. Note that $h_i^{(\beta)}$ should
be identical to $h_\ell^{(\gamma)}$ for some $\gamma\neq\beta$ and
$\ell\in[2t(\gamma)]$. Obviously, we have $s(\beta)+t(\beta
)=\mathfrak
{a}(\beta)$. We remind here the word {\emph{unpaired}} means that
$h^{(\beta)}_i$ and $h^{(\beta)}_j$ with $i\neq j$ are not identical on
$\mathcal{J}$. However, for some specific realization of $\mathbf
{j}\in
\mathcal{J}$, it is obvious that $h^{(\beta)}_i$ and $h^{(\beta)}_j$
may take the same value.
Note that since the joint cumulant is a symmetric function of the
involved variables, we can work with any specified order of these
variables. Therefore, we can write
\[
\mathbf{C}\{Z_{j_\alpha}\}_{\alpha\in A^{(\beta)}}=\mathbf {C}(Z_{l^{(\beta)}_1},Z_{l^{(\beta)}_1},
\ldots, Z_{l^{(\beta
)}_{s(\beta
)}},Z_{l^{(\beta)}_{s(\beta)}},Z_{h^{(\beta)}_1},\ldots
,Z_{h^{(\beta
)}_{2t(\beta)}}).
\]
We mimic the notation in Section~\ref{sec4.2} to denote the underlying paired
indices and unpaired indices sequence in $\{Z_{j_\alpha}\}_{\alpha\in
A^{(\beta)}}$ by
\[
\mathbf{l}^{(\beta)}:=\bigl(l_1^{(\beta)},l_1^{(\beta)},
\ldots ,l_{s(\beta
)}^{(\beta)},l_{s(\beta)}^{(\beta)}\bigr)
\]
and
\[
\mathbf{h}^{(\beta)}:=\bigl(h_1^{(\beta)},h_2^{(\beta)},
\ldots, h_{2t(\beta
)}^{(\beta)}\bigr),
\]
respectively. In addition, for simplicity we use the notation
\[
\widehat{\mathbf{l}\mathbf{h}}^{(\beta)}:=\widehat{\mathbf
{l}^{(\beta
)}\mathbf{h}^{\beta}}
\]
and write
\[
\mathbf{C}\bigl(\widehat{\mathbf{l}\mathbf{h}}^{(\beta)}\bigr):=\mathbf{C}
\{ Z_{j_\alpha}\}_{\alpha\in A^{(\beta)}}
\]
by analogy with (\ref{108.1}). Moreover, we will use the following
notation of indices sets:
\[
\bigl\{\mathbf{l}^{(\beta)}\bigr\}:=\bigl\{l_1^{(\beta)},
\ldots,l_{s(\beta
)}^{(\beta
)}\bigr\}, \qquad\bigl\{\mathbf{h}^{(\beta)}
\bigr\}:=\bigl\{h_1^{(\beta)},\ldots, h_{2t(\beta
)}^{(\beta)}
\bigr\}.
\]
We remind here that both $\mathbf{l}^{(\beta)}$ and $\mathbf
{h}^{(\beta
)}$ indices are just $\mathbf{j}$ indices in different notation. We
will call an index in $\mathbf{h}^{(\beta)},\beta\in[m]$ as
$\mathbf
{h}$ index. And $\mathbf{l}$ index can be understood analogously.
Moreover, when we refer to the position of an $\mathbf{l}$ or $\mathbf
{h}$ index, we always mean the position of its corresponding $\mathbf
{j}$ index in the word $\mathbf{j}$. In the sequel, we will also employ
the notation
%
%e4.19 #&#
\begin{equation}
B(\mathbf{h})=\{\alpha\in B: j_\alpha\mbox{ is an }\mathbf{h}\mbox{
index}\}. \label{1218.1}
\end{equation}
Note that we can regard $h_\alpha^{(\beta)}$ and $h_\gamma^{(\beta)}$
as two different free indices in $[n]$ when we take sum over $\mathcal
{J}$. However, in Lemma~\ref{lem.106.16}, the bound on the magnitude of
$\mathbf{C}(\widehat{\mathbf{l}\mathbf{h}}^{(\beta)})$ may be different
according to whether $h_\alpha^{(\beta)}$ and $h_\gamma^{(\beta)}$ take
the same value or not. Therefore, it is necessary to decompose the
summation according to different relationships between $\mathbf{h}$
indices. For example,
%
%e4.20 #&#
\begin{equation}
\qquad \sum_{h_1}\sum_{h_2}\bigl|
\mathbf{C}(Z_{h_1},Z_{h_2})\bigr|=\sum
_{h_1=h_2}\bigl|\mathbf{C}(Z_{h_1},Z_{h_2})\bigr|+\sum
_{h_1\neq h_2}\bigl|\mathbf {C}(Z_{h_1},Z_{h_2})\bigr|.
\label{1113.10}
\end{equation}
In the above example, the terms from the first summation on the
right-hand side of (\ref{1113.10}) (each term equals $1$) are quite
different from those from the second summation [each term equals
$-1/(n-1)$]. For more general $\mathbf{C}(\widehat{\mathbf{l}\mathbf
{h}})$, we will introduce the following concept of {\emph{relationship
matrix}}.

%de4.13 #&#
\begin{defi}[(Relationship matrix)] For some positive integer $N$, we
assume that $\ell_1,\ldots, \ell_N\in[n]$ and denote $\vec{\ell
}:=(\ell
_1,\ldots,\ell_N)$. Let $R_{\vec{\ell}}=(\delta_{ij})_{N,N}$ with
\[
R_{\vec{\ell}}(i,j):=\delta_{ij}=\cases{ %
 1,&\quad $\mbox{if } \ell_i=\ell_j,$
\vspace*{2pt}\cr
0,&\quad $\mbox{if } \ell_i\neq\ell_j$.}
\]
We call $R_{\vec{\ell}}$ the relationship matrix of $\vec{\ell}$.
\end{defi}

%re4.14 #&#
\begin{rem} Note that not all $0-1$ matrix can be a relationship
matrix. For example, a matrix $M$ with $M(1,2)=M(2,3)=1$ while
$M(1,3)=0$ cannot be a relationship matrix.
\end{rem}

%ex4.1 #&#
\begin{Exa} For the vector $\mathbf{j}=(1,1,1,1,2,3,3,2)$, we see the
relationship matrix of $\mathbf{j}$ is the block diagonal matrix
\begin{eqnarray*}
R_{\mathbf{j}}=\pmatrix{ 1 &1 &1 &1
\vspace*{2pt}\cr
1 &1 &1 &1
\vspace*{2pt}\cr
1 &1 &1 &1
\vspace*{2pt}\cr
1 &1 &1 &1}
\oplus\pmatrix{ 1 &0 &0 &1
\vspace*{2pt}\cr
0 &1 &1 &0
\vspace*{2pt}\cr
0 &1 &1 &0
\vspace*{2pt}\cr
1 &0 &0 &1}.
\end{eqnarray*}
\end{Exa}

%s4.5 #&#
\subsection{Proof of Proposition \texorpdfstring{\protect\ref{lem.106.15}}{4.2}}\label{sec4.5}
We recall the notation in Section~\ref{sec4.4} to write $\mathbf{C}\{
Z_{j_\alpha
}\}_{\alpha\in A^{(\beta)}}$ as
\[
\mathbf{C}\bigl(\widehat{\mathbf{l}\mathbf{h}}^{(\beta)}\bigr) =\mathbf{C}
(Z_{l_1^{(\beta)}},Z_{l_1^{(\beta)}},\ldots, Z_{l_{s(\beta)}^{(\beta)}},Z_{l_{s(\beta)}^{(\beta)}},
Z_{h_1^{(\beta
)}},Z_{h_2^{(\beta)}},\ldots, Z_{h_{2t(\beta)}^{(\beta)}} ),
\]
where $s(\beta), t(\beta)\geq0$ are nonnegative integers and
$s(\beta
)+t(\beta)=\mathfrak{a}(\beta)$. Thus, our aim is to bound the
following quantity:
%
%e4.21 #&#
\begin{equation}
\label{1014.1} \sum_{\mathbf{j}|_{B}\in\mathcal{J}(\pi_1|_{B}, n, \mathfrak
{b})}\prod
_{\beta=1}^m\bigl|\mathbf{C}\bigl(\widehat{\mathbf{l}
\mathbf{h}}^{(\beta
)}\bigr)\bigr|=\sum_{\mathbf{h}: \mathbf{j}|_{B}\in\mathcal{J}(\pi_1|_{B}, n,
\mathfrak
{b})}\prod
_{\beta=1}^m\sum_{\mathbf{l}^{(\beta)}}\bigl|
\mathbf {C}\bigl(\widehat {\mathbf{l}\mathbf{h}}^{(\beta)}\bigr)\bigr|.
\end{equation}
Here, $\sum_{\mathbf{h}: \mathbf{j}|_{B}\in\mathcal{J}(\pi
_1|_{B}, n,
\mathfrak{b})}$ represents the summation over all choices of $\mathbf
{h}$ indices along with $\mathbf{j}|_B$ running over all $\mathcal
{J}(\pi_1|_{B}, n, \mathfrak{b})$, and
%
%e4.22 #&#
\begin{equation}
\sum_{\mathbf{l}^{(\beta)}}=\sum_{l^{(\beta)}_1=1}^n
\cdots\sum_{l^{(\beta)}_{s(\beta)}=1}^n. \label{1215.1}
\end{equation}

At first, we can handle the trivial case of $m=1$ for Proposition~\ref
{lem.106.15} as follows. Observe that in this case, $A^{(1)}=B$, and
thus itself forms a block of $\pi\vee\pi_1$. Hence, $t(1)=0$ and
$s(1)=\mathfrak{a}(1)$. If $\mathfrak{a}(1)\leq2$, it is easy to see
that Proposition~\ref{lem.106.15} holds by employing (\ref{1113.20})
and Lemma~\ref{lem.1117.10}. For the case of $\mathfrak{a}(1)\geq3$,
we need to use Lemma~\ref{lem.106.16} by setting $h_1= h_2$ in (ii) and
(iii) therein. With the aid of this setting, we can get the conclusion
by noticing that except for the cases of (ii) and (iii) in Lemma~\ref
{lem.106.16} (with $h_1= h_2$ therein), the number of free indices in
any other case is at most $\mathfrak{a}(1)-2$. More specifically, we
can split the summation (\ref{1215.1}) as
\[
\sum_{\mathbf{l}^{(1)}}:=\sum^{\ast}+
\sum^{\star}+\sum^{\divideontimes},
\]
where $\sum^{\ast}$ is the summation running over the sequences
$(l^{(1)}_1,\ldots, l^{(1)}_{s(1)})\in[n]^{s(1)}$ in which all indices
are distinct from each other; $\sum^{\star}$ runs over the sequences in
which except for one pair of coincidence indices all the others are
distinct and distinct from this pair; $\sum^{\divideontimes}$ runs over
all the remaining cases. Note that the total number of the choices of
indices in $\sum^{\divideontimes}$ is $O(n^{a(1)-2})$. Then by using
(ii) and (iii) of Lemma~\ref{lem.106.16}, we can actually get the
stronger bound as $O(n^{\mathfrak{a}(1)-2})$ rather than
$O(n^{\mathfrak
{a}(1)-1})$.

Therefore, it suffices to consider the case of $m\geq2$. The proof of
this case is very complicated, so we leave it to the supplementary
material \cite{BLPZ}.

%s4.6 #&#
\subsection{A special case}\label{sec4.6}
In the sequel, we also need the following stronger bound for the
special case of $\#\{\beta: \mathfrak{a}(\beta)\geq3\}=1$ while $\#
\{
\beta:\mathfrak{a}(\beta)=2\}=0$.

%pr4.15 #&#
\begin{pro} \label{lem.1018.1}When $\#\{\beta: \mathfrak{a}(\beta
)\geq
3\}=1$ and $\#\{\beta:\mathfrak{a}(\beta)=2\}=0$ for some $B\in\pi
\vee
\pi_1$, we have
\[
\sum_{\mathbf{j}|_{B}\in\mathcal{J}(\pi_1|_{B}, n, \mathfrak{b})} \prod_{\beta=1}^{m}\bigl|
\mathbf{C}\{Z_{j_\alpha}\}_{\alpha\in A^{(\beta)}}\bigr|= O \bigl(n^{\sum_{\beta=1}^m(\mathfrak{a}(\beta)-2)\mathbf{1}_{\{
\mathfrak
{a}(\beta)\geq2\}}} \bigr).
\]
\end{pro}

The proof will be provided in the supplementary material \cite{BLPZ}.
From the above proposition, we can immediately get the following corollary.

%co4.16 #&#
\begin{cor} \label{cor.1129.10}If $\pi\in L^{\mathrm{even}}_{2k}$ such
that $\mathfrak{n}_2(\pi)=0$ and $\sum_{\gamma\geq3}\mathfrak
{n}_\gamma
(\pi)=1$, we have
\[
\sum_{\mathbf{j}\in\mathcal{J}}\bigl|\mathbf{C}_\pi(\mathbf
{j})\bigr|=O\bigl(n^{\#\pi
_1\vee\pi-1+\sum_{\gamma\geq2}(\gamma-2)\mathfrak{n}_\gamma(\pi)}\bigr).
\]
\end{cor}

%s5 #&#
\section{High order cumulants}\label{sec5}
Now with the aid of Corollaries \ref{cor.107.20}, \ref
{cor.1129.10}, and Proposition~\ref{pro.1122.1} (Proposition~3.1 of
\cite{AZ2009}) we can derive the following lemma whose proof is
provided in the supplementary material \cite{BLPZ}.

%le5.1 #&#
\begin{lem}[(High order cumulants)] \label{lem.1030.1} When $n\to
\infty$,
we have
\[
\mathbf{C} \bigl(\operatorname{tr} S^{k_1}_n,\ldots,
\operatorname{tr} S^{k_r}_n\bigr)\rightarrow 0\qquad
\mbox{for all } r\geq3.
\]
\end{lem}

%s6 #&#
\section{Mean and covariance functions}\label{sec6}
In this section, we prove (\ref{1111.2}) and (\ref{1111.3}). Before
commencing the formal proof, we introduce some necessary notation and
results on the population covariance matrix $\mathbb{E}\mathbf
{Z}^*\mathbf{Z}$ at first.
%s6.1 #&#
\subsection{On the population matrix $\mathbb{E}\mathbf{Z}^*\mathbf{Z}$}\label{sec6.1}

Let
\[
T:=T_{n,n} =\mathbb{E}\mathbf{Z}^*\mathbf{Z}=\pmatrix{ 1 &-\displaystyle\frac{1}{n-1} &\cdots&-\displaystyle\frac{1}{n-1}
\vspace*{2pt}\cr
-\displaystyle\frac{1}{n-1} &1 &\cdots&-\displaystyle\frac{1}{n-1}
\vspace*{2pt}\cr
\vdots&\vdots&\ddots&\vdots
\vspace*{2pt}\cr
-\displaystyle\frac{1}{n-1} &-\displaystyle\frac{1}{n-1} &\cdots&1 }_{n\times n}.
\]
Note that $T$ has one multiple eigenvalue $\frac{n}{n-1}$ with
multiplicity $n-1$ and one eigenvalue $0$. Roughly speaking, our aim is
to find some reference sample covariance matrix of the form $\frac
{1}{p}\Xi T\Xi^*$, where $\Xi:=(\xi_{ij})_{p,n}$ is a random matrix
with i.i.d. mean zero variance one entries, and then compare the mean
and covariance functions of the spectral statistics of $S_n$ to those
of $\frac{1}{p}\Xi T\Xi^*$. For the latter, we can use the existing
results from \cite{BS2004} and \cite{PZ2008} to obtain the explicit
formulae of the mean and covariance functions. To this end, we need to
present some notions and properties on $T$ at first. Now we denote the
empirical spectral distribution of $T$ by
\[
H_n(x):=\frac{n-1}{n}\mathbf{1}\biggl(x\geq\frac{n}{n-1}
\biggr)+\frac
{1}{n}\mathbf {1}(x\geq0).
\]
Let $c_n=n/p$ and $m_n(z):\mathbb{C}^+\to\mathbb{C}^+$ satisfy
\begin{eqnarray*}
m_n(z)&=&\int\frac{1}{t(1-c_n-c_nzm_n(z))-z}\,dH_n(t)
\\
&=&\frac{n-1}{n} \biggl[\frac{n}{n-1}\bigl(1-c_n-c_n
z m_n(z)\bigr)-z \biggr]^{-1}-\frac{1}{nz}.
\end{eqnarray*}
Regarding $m_n(z)$ as a Stieltjes transform, we can denote
$F_{c_n,H_n}$ as its corresponding distribution function.
Moreover, we denote that $H(x)=\mathbf{1}_{\{x\geq1\}}$. Define
$m(z):\mathbb{C}^+\to\mathbb{C}^+$ by the equation
$
m(z)= [(1-c-c z m(z))-z ]^{-1}
$
and set
$
\underline{m}(z)=-\frac{1-c}{z}+cm(z)$.
In order to use the results in \cite{PZ2008} (Theorem~1.4 therein), we
need to verify the following lemma on $T$, which will be proved in the
supplementary material \cite{BLPZ}.

%le6.1 #&#
\begin{lem} \label{lem.1202.2}Under the above notation, for any fixed
$z,z_1,z_2\in\mathbb{C}^+$ we have
\begin{eqnarray*}
&&\frac{1}{n}\sum_{i=1}^n
\mathbf{e}_i^*T^{1/2}\bigl(\underline {m}(z_1)T+I
\bigr)^{-1}T^{1/2}\mathbf{e}_i\mathbf{e}_i^*T^{1/2}
\bigl(\underline {m}(z_2)T+I\bigr)^{-1}T^{1/2}
\mathbf{e}_i
\\
&&\qquad\rightarrow\bigl(\underline{m}(z_1)+1\bigr)^{-1}\bigl(
\underline{m}(z_2)+1\bigr)^{-1}
\end{eqnarray*}
and
\begin{eqnarray*}
&&\frac{1}{n}\sum_{i=1}^n
\mathbf{e}_i^*T^{1/2}\bigl(\underline {m}(z)T+I
\bigr)^{-1}T^{1/2}\mathbf{e}_i\mathbf{e}_i^*T^{1/2}
\bigl(\underline {m}(z)T+I\bigr)^{-2}T^{1/2}
\mathbf{e}_i
\\
&&\qquad\rightarrow\bigl(\underline{m}(z)+1\bigr)^{-3}
\end{eqnarray*}
when $n\rightarrow\infty$. Here, $\mathbf{e}_i$ is the $n\times1$
vector with a $1$ in the $i$th coordinate and 0's elsewhere.
\end{lem}

%s6.2 #&#
\subsection{Mean function}\label{sec6.2}
At first we will pursue an argument analogous to that in Section~\ref
{sec5} to
discard the negligible terms by using Corollaries \ref{cor.107.20},
\ref{cor.1129.10}, and Proposition~\ref{pro.1122.1}. And then
we will evaluate the main terms by a {\emph{two-step comparison
strategy}} whose meaning will be clear later.

We commence with the negligible terms.
Note that for the mean function, we have $r=\#\pi_0\vee\pi_1=1$. Hence,
we can write
%
%e6.1 #&#
\begin{equation}
\mathbb{E}\operatorname{tr} S^k_n=\sum
_{\pi\in L^{\operatorname
{even}}_{2k}}p^{-k+\#\pi
_0\vee\pi}\sum_{\mathbf{j}\in\mathcal{J}}
\mathbf{C}_\pi(\mathbf{j}). \label{1023.2}
\end{equation}
Now we use bound (S.54), (S.58) and (S.60) in the
supplementary material \cite{BLPZ}. Note that when $\sum_{\gamma\geq
2}\mathfrak{n}_\gamma(\pi)\geq2$, we can easily get from (S.61) that
\[
p^{-k+\#\pi_0\vee\pi}\sum_{\mathbf{j}\in\mathcal{J}}\mathbf
{C}_\pi (\mathbf{j})=O\bigl(n^{-1}\bigr).
\]
Now we consider the case of $\sum_{\gamma\geq2}\mathfrak{n}_\gamma
(\pi)=1$.
Note that, if $\mathfrak{n}_\gamma(\pi)=1$ for any $\gamma\geq3$, we
can use the improved bound in Corollary~\ref{cor.1129.10} to obtain that
\[
p^{-k+\#\pi_0\vee\pi}\sum_{\mathbf{j}\in\mathcal{J}}\mathbf
{C}_\pi (\mathbf{j})=O\bigl(n^{-1}\bigr).
\]
However, the case of $\mathfrak{n}_2(\pi)=1$, $\mathfrak{n}_\gamma
=0,\gamma\geq3$ does have a nonnegligible contribution to the
expectation. Obviously, in this case, by (S.54) and (S.58) we see that only those terms with $\pi$ satisfying
%
%e6.2 #&#
\begin{equation}
\#\pi\vee\pi_0+\#\pi\vee\pi_1=k \label{1122.5}
\end{equation}
have $O(1)$ contribution to the total sum. Now we recall the notation
$L_{2k}^2$ and $L_{2k}^4$ defined in Section~\ref{sec3}.
We can write
%
%e6.3 #&#
\begin{eqnarray} \label{1101.1}
\mathbb{E}\operatorname{tr} S^k_n&=&\sum
_{\pi\in L_{2k}^2}p^{-k+\#\pi
_0\vee\pi
}\sum_{\mathbf{j}\in\mathcal{J}}
\mathbf{C}_\pi(\mathbf{j})
\nonumber
\\[-8pt]
\\[-8pt]
\nonumber
&&{}+\sum_{\pi\in
L_{2k}^{4}}p^{-k+\#\pi_0\vee\pi}
\sum_{\mathbf{j}\in\mathcal
{J}}\mathbf {C}_\pi(
\mathbf{j})+o(1).
\end{eqnarray}

We will not estimate the right-hand side of (\ref{1101.1}) by
bare-handed calculation and enumeration. Instead, we will adopt a
comparison approach. To this end, we need to recall some existing
results on the sample covariance matrices. At first, we define a
reference matrix. Let $\xi,\xi_j,j=1,\ldots,n$ be i.i.d. symmetric
random variables with common mean zero, variance $1$ and fourth moment
$\nu_4$. Let $V=(\xi_1,\ldots,\xi_n)$. Moreover, for any fixed positive
integer $\ell$, we assume $\mathbb{E}|\xi|^{\ell}\leq C_\ell$ for some
positive constant $C_\ell$. Then we set
$\mathbf{Y}=(Y_1,\ldots,Y_n):=VT^{1/2}$ and let $V_i,i=1,\ldots,p$ be
i.i.d. copies of $V$.
Now let $\Xi$ be the $p\times n$ matrix with $V_i$ as its $i$th row
and let
\[
S_n(\xi)=\frac{1}{p}\Xi T\Xi^*.
\]
Actually, if we take an analogous discussion on $S_n(\xi)$ as that on
$S_n$ in the last sections, it is not difficult to see that there exists
%
%e6.4 #&#
\begin{eqnarray}\label{1122.7}
\mathbb{E}\operatorname{tr} S_n^k(\xi)&=&\sum
_{\pi\in
L_{2k}^2}p^{-k+\#\pi_0\vee
\pi}\sum_{\mathbf{j}\in\mathcal{J}}
\mathbf{C}_\pi(\mathbf {j},\xi)
\nonumber
\\[-8pt]
\\[-8pt]
\nonumber
&&{}+\sum_{\pi\in L_{2k}^{4}}p^{-k+\#\pi_0\vee\pi}
\sum_{\mathbf{j}\in
\mathcal
{J}}\mathbf{C}_\pi(\mathbf{j},
\xi)+o(1),
\end{eqnarray}
where $\mathbf{C}_\pi(\mathbf{j},\xi)$ represents the quantity obtained
by replacing $Z_i$ by $Y_i$ in the definition of $\mathbf{C}_\pi
(\mathbf{j})$.
Particularly, when $\xi$ is Gaussian, we will write $\mathbb
{E}\operatorname
{tr} S_n^k(\xi)$ and $\mathbf{C}_\pi(\mathbf{j},\xi)$ as $\mathbb
{E}\operatorname{tr} S_n^k(g)$ and $\mathbf{C}_\pi(\mathbf{j},g)$,
respectively. Actually, since $\mathbf{Y}$ is just a linear transform
of i.i.d. random sequence, the verification of Lemmas \ref{lem.1011.2},
\ref{lem.1117.10} and \ref{lem.106.16} for the vector $\mathbf{Y}$ is
much easier than that for $\mathbf{Z}$. The proofs of these technical
results for $\mathbf{Y}$ are easily manipulated by invoking the
properties P1--P3 of joint cumulants stated in Section~\ref{sec3}.
However, a
more direct way to derive (\ref{1122.7}) is to check the property of
joint cumulant summability for $\mathbf{Y}$. We sketch it as follows.
At first, it is elementary check that the diagonal entries of $T^{1/2}$
are $t_d:=\sqrt{(n-1)/n}$ and the off-diagonal entries are
$t_o:=-\sqrt
{1/n(n-1)}$. Then $Y_i=t_d\xi_i+t_o\sum_{\ell\neq i} \xi_\ell$ by
definition. Now let $r$ be a fixed positive integer. Suppose that in
the collection of indices $\{j_1,\ldots,j_r\}\in[n]^{r}$, there are
$r_1\geq0$ indices taking value of $1$, and the remaining $r-r_1$
indices totally take $\alpha-1$ distinct values with multiplicities
$r_\beta\geq1, \beta=2,\ldots,\alpha$ such that $\sum_{\beta
=1}^\alpha
r_\beta=r$. Then by symmetry of $\mathbf{Y}$ and the properties P1--P3
of joint cumulant we have
\begin{eqnarray*}
&&\mathbf{C}(Y_1,Y_{j_1},\ldots, Y_{j_r})\\
&&\qquad=
\mathbf{C}(\underbrace {Y_1,\ldots,Y_1}_{r_1+1},
\underbrace{Y_2,\ldots,Y_2}_{r_2},\ldots,
\underbrace{Y_\alpha,\ldots,Y_\alpha}_{r_\alpha})
\\
&&\qquad=t_d^{1+r_1}t_o^{\sum_{\beta\neq1} r_\beta}\mathbf {C}(
\underbrace{\xi _1,\ldots,\xi_1}_{r+1})+\sum
_{\gamma=2}^\alpha t_d^{r_\gamma
}t_o^{1+\sum
_{\beta\neq\gamma}r_\beta}
\mathbf{C}(\underbrace{\xi_\gamma ,\ldots,\xi _\gamma}_{r+1}).
\end{eqnarray*}
Obviously, the quantities
$|\mathbf{C}(\underbrace{\xi_\gamma,\ldots,\xi_\gamma
}_{r+1})|,\gamma
=1,\ldots, \alpha$ are all the same and can be bounded by some positive
constant from above by invoking the assumption that $\mathbb{E}|\xi
|^{\ell}\leq C_\ell$ and the formula (\ref{106.2}). Moreover, we
observe that $\sum_{\beta\neq1}r_\beta\geq\alpha-1$ and $1+\sum_{\beta
\neq\gamma}r_\beta\geq\alpha-1$ for all $\gamma=2,\ldots,\alpha
$. In
addition, we have $t_d=O(1)$, $t_o=O(n^{-1})$. Therefore,
\[
\bigl|\mathbf{C}(Y_1,Y_{j_1},\ldots, Y_{j_r})\bigr|=O
\bigl(n^{-\alpha+1}\bigr).
\]
Observe that $\alpha-1$ is the number of distinct values except for $1$
in the collection $\{j_1,\ldots,j_r\}$. Consequently, we have that
\[
\sum_{j_1=1}^n\cdots\sum
_{j_r=1}^n\bigl|\mathbf{C}(Y_1,Y_{j_1},
\ldots, Y_{j_r})\bigr|=O(1),
\]
which implies that the joint cumulant summability holds for $\mathbf
{Y}$. As explained above, we can get that the stronger bound
\[
\sum_{\mathbf{j}\in\mathcal{J}}\bigl|\mathbf{C}_\pi(\mathbf{j},
\xi )\bigr|= O\bigl(n^{\#
\pi_1\vee\pi}\bigr)
\]
holds by using the result in \cite{AZ2009}. With the aid of this
stronger bound, we can derive (\ref{1122.7}) by a routine discussion
as before.

In addition, obviously, we have
%
%e6.5 #&#
\begin{equation}
\sum_{\pi\in L_{2k}^2}p^{-k+\#\pi_0\vee\pi}\sum
_{\mathbf{j}\in
\mathcal
{J}}\mathbf{C}_\pi(\mathbf{j},\xi)=\sum
_{\pi\in L_{2k}^2}p^{-k+\#
\pi
_0\vee\pi}\sum_{\mathbf{j}\in\mathcal{J}}
\mathbf{C}_\pi(\mathbf{j}) \label{1215.124}
\end{equation}
since this term only depends on the covariance structure $T$ which is
shared by $\mathbf{Z}$ and $\mathbf{Y}$.

For the second term on the right-hand side of (\ref{1122.7}), by
Corollary~\ref{cor.107.20} and (\ref{1122.5}) we see that it has a
contribution of $O(1)$ at most. Hence, it suffices to estimate its
leading term. To this end, we use (\ref{1024.3}) with $\mathbf{Z}$
replaced by $\mathbf{Y}$. Now we consider the sum over $\mathbf
{j}|_{B_i}$ with the block $B_i$ containing the unique $4$-element
block $A_i^{(\gamma)}\in\pi$. Obviously the sums over the indices with
positions in the other blocks of $\pi\vee\pi_1$ are all in the case of
(\ref{1129.5}). Without loss of generality, we can fix $i$ and $\gamma$
in the following argument. Recall the notation of paired index and
unpaired index. Observe that $\mathbf{C}\{Y_{j_\alpha}\}_{\alpha\in
A_i^{(\gamma)}}$ may be in one of the following forms. When $m_i=1$,
$\mathbf{C}\{Y_{j_\alpha}\}_{\alpha\in A_i^{(\gamma)}}$ must be in the
form of $\mathbf{C}(Y_{l^{(\gamma)}_1},Y_{l^{(\gamma
)}_1},Y_{l^{(\gamma
)}_2},Y_{l^{(\gamma)}_2})$ (case 1). When $m_i\geq2$, it may be in the
form of $\mathbf{C}(Y_{l^{(\gamma)}_1},Y_{l^{(\gamma
)}_1},Y_{h^{(\gamma
)}_1},Y_{h^{(\gamma)}_2})$ (case 2) or $\mathbf{C}(Y_{h^{(\gamma
)}_1},Y_{h^{(\gamma)}_2},Y_{h^{(\gamma)}_3},Y_{h^{(\gamma)}_4})$ (case
3). Then we have the following lemma which will be proved in
supplementary material \cite{BLPZ}.

%le6.2 #&#
\begin{lem} \label{lemma5.2} In any of the above three cases, for given
$\pi|_{B_i}$ and $\pi_1|_{B_i}$, there is a unique $\sigma^{(\gamma
)}\in L_4^2$ such that
%
%e6.6 #&#
\begin{equation}
\prod_{\beta\neq\gamma}\mathbf{C}\{Y_{j_\alpha}\}
_{\alpha\in A_i^{(\beta)}}\cdot\mathbf{C}_{\sigma^{(\gamma)}}\{ Y_{j_\alpha}
\}_{\alpha\in A_i^{(\gamma)}} \label{1216.1}
\end{equation}
is a product of two cycles. The other two perfect matchings in $L_4^2$
will drive the above product to be only one cycle.
\end{lem}

Note that by Proposition~\ref{lem.106.15}, we have
%
%e6.7 #&#
\begin{equation}
\sum_{\mathbf{j}|_{B_i}\in\mathcal{J}(\pi_1|_{B_i}, n, \mathfrak
{b}_i)}\prod_{\beta=1}^{m_i}\bigl|
\mathbf{C}\{Y_{j_\alpha}\}_{\alpha\in
A_i^{(\beta)}}\bigr|\leq O(n). \label{1122.6}
\end{equation}
Our aim is to get the explicit $O(n)$ term of (\ref{1122.6}). To this
end, we recall the discussions in Section~\ref{sec4}. In any case of
$\mathbf
{C}\{Y_{j_\alpha}\}_{\alpha\in A_i^{(\gamma)}}$, we only need to
consider those $\mathbf{j}$ such that $(j_\alpha)_{\alpha\in
A_i^{(\gamma)}}$ is $\sigma^{(\gamma)}$-measurable. Since we can see
that all the other $\mathbf{j}|_{B_i}$ can only make an $O(1)$
contribution totally to the left-hand side of (\ref{1122.6}) by using
Lemma~\ref{lem.1117.10} and the discussion on summations of
in-homogeneous cycles in Section~\ref{sec4.3}. Moreover, by the
discussions in
Section~\ref{sec4}, we know that in any of the aforementioned three
cases, the
main contribution to the summation comes from the terms which can be
decomposed into homogeneous cycles. In these terms, each $2$-element
cumulant is equal to $1$. Hence, we have
%
%e6.8 #&#
\begin{eqnarray}
\label{1122.9} &&\sum_{\mathbf{j}|_{B_i}\in\mathcal{J}(\pi_1|_{B_i}, n, \mathfrak
{b}_i)}\mathbf{C}
\{Y_{j_\alpha}\}_{\alpha\in A_i^{(1)}}\nonumber\\
 &&\qquad=\sum_{\alpha_1,\alpha_2=1}^n
\mathbf{C}(Y_{\alpha_1},Y_{\alpha
_1},Y_{\alpha_2},Y_{\alpha_2})+O(1)
\\
&&\qquad =n\mathbf{C}(Y_1,Y_1,Y_1,Y_1)+n(n-1)
\mathbf {C}(Y_1,Y_1,Y_2,Y_2)+O(1).\nonumber
\end{eqnarray}
Note that (\ref{1122.9}) also holds if we replace $\mathbf{Y}$
variables by corresponding $\mathbf{Z}$ variables.

Now for the Gaussian case we claim
%
%e6.9 #&#
\begin{equation}
\sum_{\pi\in L_{2k}^{4}}p^{-k+\#\pi_0\vee\pi}\sum
_{\mathbf{j}\in
\mathcal
{J}}\mathbf{C}_\pi(\mathbf{j},g)=o(1).
\label{1101.2}
\end{equation}
To see (\ref{1101.2}), it suffices to show that for any fixed $\pi\in
L^4_{2k}$, (\ref{1122.6}) can be strengthened to be
%
%e6.10 #&#
\begin{equation}
\sum_{\mathbf{j}|_{B_i}\in\mathcal{J}(\pi_1|_{B_i}, n, \mathfrak
{b}_i)}\prod_{\beta=1}^{m_i}\bigl|
\mathbf{C}\{Y_{j_\alpha}\}_{\alpha\in
A_i^{(\beta)}}\bigr|\leq O(1) \label{1215.123}
\end{equation}
when $\xi$ is Gaussian. By (\ref{1122.9}), it suffices to evaluate the
quantity $\mathbf{C}(Y_{\alpha_1},Y_{\alpha_1},\break Y_{\alpha
_2},  Y_{\alpha_2})$.
Note that when $\alpha_1=\alpha_2$,
%
%e6.11 #&#
\begin{equation}
\mathbf{C}(Y_{\alpha_1},Y_{\alpha_1},Y_{\alpha_2},Y_{\alpha
_2})=
\mathbf {C}(Y_{1},Y_{1},Y_{1},Y_{1})=
\nu_4-3=0 \label{1101.3}
\end{equation}
since $Y_1$ is Gaussian. If $\alpha_1\neq\alpha_2$, it is not
difficult to get that
%
%e6.12 #&#
\begin{equation}
\mathbf{C}(Y_{\alpha_1},Y_{\alpha_1},Y_{\alpha_2},Y_{\alpha
_2})=
\mathbf {C}(Y_{1},Y_{1},Y_{2},Y_{2})=O
\bigl(n^{-2}\bigr) \label{1101.4}
\end{equation}
by the definition of $Y_i$ and propositions P1--P3 of joint cumulant.
Thus, (\ref{1215.123}) holds, which directly implies that
\[
\sum_{\mathbf{j}\in\mathcal{J}}\mathbf{C}_\pi(\mathbf
{j},g)=O\bigl(n^{\#\pi
_1\vee\pi-1}\bigr),\qquad \pi\in L^4_{2k}.
\]
It further yields (\ref{1101.2}) by combining (S.58) and the
elementary fact that $\#\pi=k-1$ for $\pi\in L_{2k}^4$.
Inserting (\ref{1215.124}) and (\ref{1101.2}) into (\ref{1122.7}),
we obtain
\[
\sum_{\pi\in L_{2k}^2}p^{-k+\#\pi_0\vee\pi}\sum
_{\mathbf{j}\in
\mathcal
{J}}\mathbf{C}_\pi(\mathbf{j})=\mathbb{E}
\operatorname{tr} S^k_n(g)+o(1).
\]
Therefore, for the first term on the right-hand side of (\ref{1101.1}),
it suffices to estimate $\mathbb{E}\operatorname{tr} S^k_n(g)$. For
the latter,
we can use the result from \cite{BS2004} or \cite{PZ2008} directly to
write down
%
%e6.13 #&#
\begin{eqnarray}\qquad
\label{1112.1} &&\sum_{\pi\in L_{2k}^2}p^{-k+\#\pi_0\vee\pi}\sum
_{\mathbf{j}\in
\mathcal
{J}}\mathbf{C}_\pi(\mathbf{j})
\nonumber
\\[-8pt]
\\[-8pt]
\nonumber
&&\qquad=n\int x^k\,d F_{c_n,H_n}(x)-\frac{1}{2\pi i}\int
_{\mathcal{C}} \frac
{cz^k\underline{m}^3(z)/(1+\underline{m}(z))^3}{ [1-c\underline
{m}^2(z)/(1+\underline{m}(z))^2 ]}\,dz+o(1),
\end{eqnarray}
where the contour $\mathcal{C}$ is taken to enclose the interval
$[(1-\sqrt{c})^2,(1+\sqrt{c})^2]$ as interior. See Theorem~1.4 of
\cite
{PZ2008}, for instance.
However, here we can further simplify (\ref{1112.1}) by the property of
orthogonal invariance of standard Gaussian vectors. Note that when $\xi
$ is Gaussian, we have
\[
\frac{1}{p}\Xi T \Xi^*\stackrel{d}=\frac{1}{p}
\frac{n}{n-1}GG^*,
\]
where $G:=(g_{i,j})_{p,n-1}$ with i.i.d. $N(0,1)$ elements. Now let
$\tilde{c}_n=\frac{n-1}{p}$ and $F_{\tilde{c}_n,MP}$ be
Marchenko--Pastur law (MP law) with parameter $\tilde{c}_n$. Then by
Theorem~1.4 of~\cite{PZ2008} and Lemma~\ref{lem.1202.2} we can rewrite
(\ref{1112.1}) as
\begin{eqnarray*}
\sum_{\pi\in L_{2k}^2}p^{-k+\#\pi_0\vee\pi}\sum
_{\mathbf{j}\in
\mathcal
{J}}\mathbf{C}_\pi(\mathbf{j})&=&\frac{n^k}{(n-1)^{k-1}}
\int x^k\,d F_{\tilde{c}_n,MP}(x)
\nonumber
\\[-8pt]
\\[-8pt]
\nonumber
&&{}-\frac{1}{2\pi i}\int_{\mathcal{C}} \frac{cz^k\underline
{m}^3(z)/(1+\underline{m}(z))^3}{ [1-c\underline
{m}^2(z)/(1+\underline{m}(z))^2 ]}\,dz+o(1).
\end{eqnarray*}
Note that by (9.8.14) of \cite{BS2004}, we see that
\begin{eqnarray*}
&&-\frac{1}{2\pi i}\int_{\mathcal{C}} \frac{cz^k\underline
{m}^3/(1+\underline{m}(z))^3}{ [1-c\underline
{m}^2(z)/(1+\underline
{m}(z))^2 ]}\\
&&\qquad=
\frac{1}{4} \bigl[(1-\sqrt{c})^{2k}+(1+\sqrt {c})^{2k}
\bigr]-\frac{1}{2}\sum_{j=0}^k
\pmatrix{k
\cr
j}^2c^j.
\end{eqnarray*}
Moreover, by the formula of moments of MP law (see Section~3.1.1 of
\cite{BS2004}, e.g.) one can also get that
\[
\frac{n^k}{(n-1)^{k-1}}\int x^k\,d F_{\tilde{c}_n,MP}(x)=\frac
{n^k}{(n-1)^{k-1}}
\sum_{j=0}^{k-1}\frac{1}{j+1}\pmatrix{k
\cr
j}\pmatrix {k-1
\cr
j} \biggl(\frac{n-1}{p} \biggr)^{j}.
\]

Now we come to estimate the second term on the right-hand side of (\ref
{1101.1}). Now we choose $\xi$ satisfying $\nu_4\neq3$. Note that
(\ref
{1101.3}) is not valid now. However, (\ref{1101.4}) still holds. Then
in this case,
\[
(\ref{1122.9})=(\nu_4-3) n+O(1).
\]
Note that according to (\ref{1024.3}), except for this $B_i$ which
containing the unique $4$-element block of $\pi$, the summation over
the indices with positions in $[2k]\setminus B_i$ only depends on the
covariance structure since $\pi\in L_{2k}^4$. Now for
$C(Z_{l_1},Z_{l_1},Z_{l_2},Z_{l_2})$, we have
\[
C(Z_{l_1},Z_{l_1},Z_{l_2},Z_{l_2})=
\mathbf{E}Z_1^4-3+o(1), \qquad l_1=l_2
\]
and
\[
C(Z_{l_1},Z_{l_1},Z_{l_2},Z_{l_2})=
\frac{1}{n}\bigl(1-\mathbf {E}Z_1^4\bigr)+O
\bigl(n^{-2}\bigr),\qquad l_1\neq l_2,
\]
which can be checked easily by the distribution of $\mathbf{Z}$. Then
we have
\[
(\ref{1122.9})_{[\mathbf{Y}\to\mathbf{Z}]}=-2n+O(1),
\]
where $(\ref{1122.9})_{[\mathbf{Y}\to\mathbf{Z}]}$ represents the
quantity obtained by replacing $\mathbf{Y}$ by $\mathbf{Z}$ in (\ref
{1122.9}).
And all the other factors in (\ref{1024.3}) are the same as those of
$\sum_{\mathbf{j}\in\mathcal{J}}\mathbf{C}_\pi(\mathbf{j},\xi)$ since
they only depends on the covariance matrix $T$. That means
\[
\sum_{\pi\in L_{2k}^{\operatorname{4}}}p^{-k+\#\pi_0\vee\pi}\sum
_{\mathbf{j}\in
\mathcal{J}}\mathbf{C}_\pi(\mathbf{j})=-
\frac{2}{\nu_4-3}\sum_{\pi\in
L_{2k}^{\operatorname{4}}}p^{-k+\#\pi_0\vee\pi}\sum
_{\mathbf
{j}\in\mathcal
{J}}\mathbf{C}_\pi(\mathbf{j},
\xi)+o(1).
\]
By using \cite{PZ2008} again (see Theorem~1.4 therein), we can get that
\[
\sum_{\pi\in L_{2k}^{\operatorname{4}}}p^{-k+\#\pi_0\vee\pi}\sum
_{\mathbf{j}\in
\mathcal{J}}\mathbf{C}_\pi(\mathbf{j},\xi)=
\frac{1}{\pi i}\int_{\mathcal
{C}}\frac{cz^k\underline{m}^3(z)(\underline{m}(z)+1)^{-3}}{
1-c\underline{m}^2(z)/(1+\underline{m}(z))^2}\,dz+o(1).
\]
By (1.23) of \cite{PZ2008}, we obtain
\begin{eqnarray*}
\frac{1}{\pi i}\int_{\mathcal{C}}\frac{cz^k\underline
{m}^3(z)(\underline
{m}(z)+1)^{-3}}{
1-c\underline{m}^2(z)/(1+\underline{m}(z))^2}\,dz&=&2c^{1+k}
\sum_{j=0}^k\pmatrix{k
\cr
j} \biggl(
\frac{1-c}{c} \biggr)^j\pmatrix {2k-j
\cr
k-1}
\\
&&{}- 2c^{1+k}\sum_{j=0}^k
\pmatrix{r
\cr
j} \biggl(\frac{1-c}{c} \biggr)^j\pmatrix {2k+1-j
\cr
k-1}.
\end{eqnarray*}

In summary, we use the Gaussian matrix as the reference matrix to
obtain the value of the summation over $\pi\in L_{2k}^2$ and then use
the general matrix with $\nu_4\neq3$ as the reference one to obtain
the value of the summation over $\pi\in L_{2k}^4$. We call such a
comparison strategy as a {\emph{two-step comparison strategy}}.
%s6.3 #&#
\subsection{Covariance function}\label{sec6.3}
Now we estimate the covariance function. Again we start with the formula
\[
\mathbf{C}\bigl(\operatorname{tr} S_n^{k_1},\operatorname{tr} S_n^{k_2}\bigr)
=\mathop{\sum_{\pi\in L_{2k}^{\operatorname{even}}}}_{
\operatorname{s.t.} \#\pi_0\vee
\pi_1\vee\pi=1}p^{-k+\#\pi_0\vee\pi}\sum_{\mathbf{j}\in
\mathcal
{J}}\mathbf{C}_\pi(\mathbf{j}).
\]
Similar to the discussion in the last subsection, by using Corollaries
\ref{cor.107.20} and \ref{cor.1129.10} and (S.55) again we can
see that is suffices to evaluate the contribution of the summation over
the partitions $\pi$ satisfying $\mathfrak{n}_2(\pi)=0\mbox{ or }1$ and
$\mathfrak{n}_\gamma(\pi)=0$ for all $\gamma\geq3$. Moreover, by
(S.54) and (S.55) it is easy to see that
$
\mathbf{C}(\operatorname{tr} S_n^{k_1},\operatorname{tr} S_n^{k_2})=O(1)
$
since $\#\pi_0\vee\pi+\#\pi_1\vee\pi\leq k$ when $r=2$. Now let
\begin{eqnarray*}
\tilde{L}_{2k}^2&:=&\bigl\{\pi\in L_{2k}^2:
\#\pi\vee\pi_1\vee\pi_0=1\bigr\},\\
\tilde{L}_{2k}^4&:=&\bigl\{\pi\in L_{2k}^4:
\#\pi\vee\pi_1\vee\pi_0=1\bigr\}.
\end{eqnarray*}

For the explicit evaluation, we adopt the aforementioned {\emph
{two-step comparison strategy}} again. We split the summation into the
summations over $\tilde{L}_{2k}^2$ partitions and $\tilde{L}_{2k}^4$
partitions. For the first part, we compare it with that of the Gaussian
case. And for the second part, we compare it with the case of $\nu
_4\neq3$. Then it is analogous to use Theorem~1.4 of \cite{PZ2008} and
Lemma~\ref{lem.1202.2} to obtain that
\begin{eqnarray*}
&&\mathbf{C}\bigl(\operatorname{tr} S_n^{k_1},
\operatorname{tr} S_n^{k_2}\bigr)\\
&&\qquad=-\frac{1}{2\pi
^2}\int
_{\mathcal{C}_1}\int_{\mathcal{C}_2}\frac
{z^{k_1}_1z^{k_2}_2}{(\underline{m}(z_1)-\underline
{m}(z_2))^2}
\underline{m}'(z_1)\underline{m}'(z_2)\,dz_1\,dz_2
\\
&&\qquad\quad{} +\frac{c}{2\pi^2}\int_{\mathcal{C}_1}\int_{\mathcal{C}_2}
z_1^{k_1}z_2^{k_2}\frac{d^2}{dz_1\,dz_2}
\biggl(\frac{\underline
{m}(z_1)\underline{m}(z_2)}{(\underline{m}(z_1)+1)(\underline
{m}(z_2)+1)} \biggr)\,dz_1\,dz_2+o(1),
\end{eqnarray*}
where the contours $\mathcal{C}_1$ and $\mathcal{C}_2$ are disjoint and
both enclose the interval $[(1-\sqrt{c})^2,(1+\sqrt{c})^2]$ as interior.
Now by (9.8.15) of \cite{BS2004} and (1.24) of \cite{PZ2008}, we have
\begin{eqnarray*}
&&-\frac{1}{2\pi^2}\int_{\mathcal{C}_1}\int_{\mathcal{C}_2}
\frac
{z^{k_1}_1z^{k_2}_2}{(\underline{m}(z_1)-\underline
{m}(z_2))^2}\underline{m}'(z_1)
\underline{m}'(z_2)\,dz_1\,dz_2
\nonumber
\\
&&\qquad=2c^{k_1+k_2}\sum_{j_1=0}^{k_1-1}\sum
_{j_2=0}^{k_2}\pmatrix {k_1
\cr
j_1}\pmatrix{k_2
\cr
j_2} \biggl(
\frac{1-c}{c} \biggr)^{j_1+j_2}
\nonumber
\\
&&\qquad\quad{}\times\sum_{l=1}^{k_1-j_1}l
\pmatrix{2k_1-1-(j_1+l)
\cr
k_1-1}\pmatrix
{2k_2-1-j_2+l
\cr
k_2-1}
\end{eqnarray*}
and
\begin{eqnarray*}
\hspace*{-4pt}&&\frac{c}{2\pi^2}\int_{\mathcal{C}_1}\int_{\mathcal{C}_2}
z_1^{k_1}z_2^{k_2}\frac{d^2}{dz_1\,dz_2}
\biggl(\frac{\underline
{m}(z_1)\underline{m}(z_2)}{(\underline{m}(z_1)+1)(\underline
{m}(z_2)+1)} \biggr)\,dz_1\,dz_2
\nonumber
\\
\hspace*{-4pt}&&\qquad=-2c^{k_1+k_2+1}\sum_{j_1=0}^{k_1}\sum
_{j_2=0}^{k_2}\pmatrix {k_1
\cr
j_1}\pmatrix{k_2
\cr
j_2} \biggl(
\frac{1-c}{c} \biggr)^{j_1+j_2}\\
&&\qquad\quad{}\times \pmatrix {2k_1-j_1
\cr
k_1-1}\pmatrix{2k_2-j_2
\cr
k_2-1}.
\end{eqnarray*}
Thus, we complete the proof of Theorem~\ref{thmm.1111.1}.

\section*{Acknowledgements}
The authors are very grateful to Professor Greg W. Anderson for his
answer to our question on Proposition~3.1 of \cite{AZ2009}. We also
thank the anonymous referees and an Associate Editor for careful
reading and valuable comments.

\begin{supplement}%[id=suppA]
\stitle{Supplement to ``Spectral statistics of large dimensional Spearman's
rank correlation matrix and its application''}
\slink[doi]{10.1214/15-AOS1353SUPP} %[doi,text={...}] - jei reikia suskaldyti doi
\sdatatype{.pdf}
\sfilename{aos1353\_supp.pdf}
\sdescription{This supplemental article \cite{BLPZ} contains the proofs of Lemmas
\ref
{lem.1011.2}, \ref{lem.105.1} \ref{lem.1117.10}, \ref{lem.106.16},
\ref
{lem.1105.2},
Propositions \ref{lem.106.15}, \ref{lem.1018.1}, Lemmas \ref
{lem.1030.1}, \ref{lem.1202.2}, \ref{lemma5.2}.}
\end{supplement}

% imsref loaded by akundreckaite, 2015-08-25 14:01:59
%

%\begin{appendix}
%\section{}
%\end{appendix}

%\section*{Acknowledgments}

%\begin{thebibliography}{99}
%\bibitem[\protect\citeauthoryear{}{}]{r1}
%\bibitem{r1}
%\end{thebibliography}

\printaddresses
\end{document}